\documentclass[11pt]{amsart}
\usepackage{amsmath}
\usepackage{amscd}
\usepackage{amssymb}
\usepackage{pdfsync}

\newtheorem{theorem}{Theorem}[section]

\newtheorem{proposition}[theorem]{Proposition}

\newtheorem*{prop}{Proposition}

\newtheorem*{lem}{Lemma}
\newtheorem*{rem}{Remark}

\theoremstyle{definition}

\newcommand{\cN}{{\mathcal N}}

\newcommand{\cP}{{\mathcal P}}

\newcommand{\cU}{{\mathcal U}}

\newcommand{\cS}{{\mathcal{S}}}

\newcommand{\bW}{{\mathbf W}}

\newcommand{\rc}{{\mathrm c}}

\newcommand{\Lg}{{\mathfrak g}}

\newcommand{\tF}{{\textbf{F}}}
\newcommand{\tk}{{\textbf{k}}}

\newcommand{\p}{\perp}

\newcommand{\la}{\langle}
\newcommand{\ra}{\rangle}

\newcommand{\beq}{\begin{equation*}}
\newcommand{\eeq}{\end{equation*}}

\begin{document}

\title[On Unipotent and Nilpotent Pieces for Classical Groups]
{On Unipotent and Nilpotent Pieces for Classical Groups}
        \author{Ting Xue}
        \address{Department of Mathematics, Northwestern University,
Evanston, IL 60208, USA}
        \email{txue@math.northwestern.edu}
\maketitle
\begin{abstract}
We show that the definition of unipotent (resp. nilpotent) pieces for classical groups given by Lusztig (resp. Lusztig and the author) coincides with the combinatorial definition using closure relations on unipotent  classes (resp. nilpotent orbits). Moreover we give a closed formula for a map from the set of unipotent  classes (resp. nilpotent orbits) in characteristic 2 to the set of unipotent classes in characteristic 0 such that the fibers are the unipotent (resp. nilpotent) pieces.
\end{abstract}

\section{Introduction}
Let $G$ be a connected reductive group of type $B$, $C$ or $D$ defined
over an algebraically closed field of characteristic exponent $p\geq 1$ and let $\Lg$ be
the Lie algebra of $G$. Denote $\cU_G$ (resp. $\cN_\Lg$) the set of
unipotent (resp. nilpotent) elements in $G$ (resp. $\Lg$).
In
\cite{Lu2,Lu3}, Lusztig defines a partition of $\cU_G$ into smooth locally closed $G$-stable pieces, called
unipotent pieces (see \cite{Lu2} for symplectic groups and
\cite{Lu3} for special orthogonal groups). In \cite{Lu4}, Lusztig proposes
another definition of unipotent pieces which  unifies the
definitions  in \cite{Lu2,Lu3}. In Appendix A of \cite{Lu4}, Lusztig and the author define an analogue partition of $\cN_\Lg$ into smooth locally closed $G$-stable pieces, called nilpotent pieces.   The unipotent or nilpotent pieces are indexed by unipotent classes in the group $G_\mathbb{C}$ over $\mathbb{C}$ of the same type as $G$, and in many
ways depend very smoothly on $p$. In particular, the number of $\tF_{p^s}$-rational points in a unipotent or nilpotent piece  is a polynomial in $p^s$ independent of $p,s$.

On the other hand, there is a natural injection map from the set of unipotent classes in $G_\mathbb{C}$ to the set of unipotent classes in $G$ (given by the Springer correspondence). Using this map and the closure relation on unipotent classes one can define a partition of $\cU_G$ into locally closed pieces, which are called MS-pieces (after Mizuno and Spaltenstein) by Lusztig \cite{Lu7}. We show in section \ref{sec-comdef} and section \ref{sec-fiber} that the MS-pieces are the same as unipotent pieces defined by Lusztig  (for symplectic groups this follows from \cite{Lu2,Lu4}). In view of the properties of unipotent pieces proved by Lusztig, this implies  that the MS-pieces are smooth and that the number of $\tF_{p^s}$-rational points in an MS-piece  is a polynomial in $p^s$ independent of $p,s$ (this is the statement 6.8 (a) of \cite{Lu7} for classical groups). We also define MS-pieces in $\cN_\Lg$ and prove analogous results for $\cN_\Lg$. In particular, we determine in
Proposition \ref{prop-psi} which unipotent classes (resp. nilpotent orbits) lie in
the same piece (for pieces in symplectic groups this follows
from \cite{Lu2,Lu4};
for unipotent pieces in special orthogonal groups another computation
using different methods is given in \cite{Lu8}).

In section \ref{sec-spp}, we define a partition of $\cU_G$ (resp. $\cN_\Lg$) into special pieces as in \cite{Lu7} (where $p=1$) and show that a special piece is a union of unipotent (resp. nilpotent) pieces (for $\cU_G$ this follows implicitly from \cite[III]{Spal}, see \cite{Lu8}). We also explain how this implies that the number of $\tF_{p^s}$-rational points in a special piece  is a polynomial in $p^s$ that depends only on the Weyl group of $G$ (this is the statement 6.9 (a) of \cite{Lu7} for classical groups). 
\section{Notations and recollections}
\subsection{Orders on the set of partitions and the set of pairs of partitions}\label{sec-partition}
Let $\cP(n)$ denote the set of all partitions
$\lambda=(\lambda_1\geq\lambda_2\geq\cdots\geq 0)$ such that
$|\lambda|:=\sum\lambda_i=n$. For $\lambda\in\cP(n)$, define $\lambda_j^*=|\{\lambda_i|\lambda_i\geq j\}|$ and $m_\lambda(j)=\lambda_{j}^*-\lambda_{j+1}^*$. For $\lambda,\mu\in\cP(n)$, we say
that $\lambda\leq\mu$ if the following equivalent conditions (a) and (a$'$) hold

\qquad (a) $\sum_{j\in[1,i]}\lambda_j\leq\sum_{j\in[1,i]}\mu_j\text{ for all }i\geq 1,$

\qquad (a$'$) $\sum_{j\in[1,i]}\lambda_j^*\geq\sum_{j\in[1,i]}\mu_j^*\text{ for all }i\geq 1.$

Let $\cP_2(n)$ denote the set of all pairs of partitions
$(\alpha,\beta)$ such that $|\alpha|+|\beta|=n$. For $(\alpha,\beta)\in\cP_2(n)$, $\alpha=(\alpha_1\geq\alpha_2\geq\cdots)$, $\beta=(\beta_1\geq\beta_2\geq\cdots)$, we set
\begin{equation}\label{eqn-ab}A_i=\sum_{j\in[1,i]}(\alpha_j+\beta_j),\quad B_i=\sum_{j\in[1,i-1]}(\alpha_j+\beta_j)+\alpha_i,\ i\geq 1.\end{equation}
For
$(\alpha,\beta),(\alpha',\beta')\in\cP_2(n)$, we say that
$(\alpha,\beta)\leq(\alpha',\beta')$ if
 $
A_i\leq A_i'$ and $B_i\leq B_i'
\text{ for all }i\geq 1.$

\subsection{Combinatorial parametrization of irreducible Weyl group characters}\label{ssec-notation}
Let $\bW$ be the Weyl group of $G$ and $\bW^\wedge$
the set of irreducible characters of $\bW$ over $\mathbb{C}$.

If $\mathbf{W}$ is of type $B_n$ (or $C_n$), $n\geq
1$, then $\mathbf{W}^\wedge$ is parametrized by ordered
pairs of partitions $(\alpha,\beta)\in\cP_2(n)$ (see \cite{Lu1}). We identify
$\mathbf{W}^\wedge$ with $\cP_2(n)$ where $(n,-)$ is the trivial character and $(-,1^n)$ is the sign character.

If $\mathbf{W}$ is of type $D_n$, $n\geq 2$, then $\mathbf{W}^\wedge$  is parametrized by unordered
pairs of partitions $\{\alpha,\beta\}$ with $|\alpha|+|\beta|=n$ where each pair $\{\alpha,\alpha\}$
corresponds to two (degenerate) elements of
$\mathbf{W}^\wedge$ (see \cite{Lu1}).
  We identify ${\mathbf{W}^\wedge}$ with the set
$\{(\alpha,\beta)\in\cP_2(n)|\beta_1\leq\alpha_1\}$ where each pair $(\alpha,\alpha)$ is counted twice.

\subsection{Combinatorial description of the images of Springer correspondence maps}

Denote $\Omega_{G}^p$  the set of unipotent classes in $\cU_G$ and
$\Omega_{\Lg}^p$ the set of nilpotent orbits in $\cN_{\Lg}$. Recall that we have injective maps (see \cite{Sp1,Lu6,X2})
\begin{equation*}\label{injmaps}
\gamma_{G}^p:\Omega_{G}^p\rightarrow \bW^\wedge,\
\gamma_{\Lg}^p:\Omega_{\Lg}^p\rightarrow \bW^\wedge
\end{equation*} which map a
class/orbit $\rc$ to the irreducible character of $\bW$ corresponding to
the pair $(\rc,1)$ under Springer correspondence.
 We denote $\Lambda_{G}^p$ (resp. $\Lambda_{\Lg}^p$) the image of the
map  $\gamma_{G}^p$ (resp.
$\gamma_{\Lg}^p$). We may write $\Omega^p$ to denote either $\Omega_{G}^p$ or $\Omega_{\Lg}^p$ (when it is easy to determine from the context) and  similar conventions apply for $\Lambda^p$, $\gamma^p$.

When $p\neq 2$, we can identify $\Omega_{\Lg}^p$ with $\Omega_{G}^p$, $\gamma_{\Lg}^p$ with $\gamma_{G}^p$, and $\Lambda_{\Lg}^p$ with $\Lambda_{G}^p$ since by Springer \cite{Sp2} there exists a $G$-equivariant isomorphism from $\cU_G$ to $\cN_\Lg$, and we can further identify $\Omega_{G}^p$ ($p\neq 2$) with  $\Lambda_{G}^1$,  $\gamma_{G}^p$ $(p\neq 2)$ with $\gamma_{G}^1$, and $\Lambda_{G}^p$ $(p\neq 2)$ with $\Lambda_{G}^1$ since the classification of unipotent classes in $G$ is the same as that in $G_\mathbb{C}$ (which depends only on the type of $G$).

Since the classification of nilpotent orbits or unipotent classes in $G$ depends only on the type of $G$,  we  assume from now on that $G$ is a symplectic group $Sp(2n)$ or a special orthogonal group $SO(N)$. We will often identify $G=Sp(2n)$ with $Sp_{\la,\ra}(V)$ (resp.  $G=SO(N)$ with $SO_Q(V)$), where $V$ is a vector space of dimension $2n$ (resp. $N$) over $\tk$ equipped with a fixed non-degenerate symplectic form $\la,\ra$ (resp. fixed
non-degenerate quadratic form $Q$), $Sp_{\la,\ra}(V)=\{g\in GL(V)|\la gv, gv'\ra=\la v,v'\ra,\ \forall\ v,v'\in V\}$ (resp. $SO_Q(V)$ is the identity component of $O_Q(V)=\{g\in
GL(V)|Q(gv)=v,\forall\ v\in V\}$). Thus  $\Lg=\mathfrak{sp}(2n)=\mathfrak{sp}_{\la,\ra}(V)=\{x\in\mathfrak{gl}(V)|\la xv,v'\ra+\la v,x v'\ra=0,\ \forall\ v,v'\in V\}$ (resp. $\Lg=\mathfrak{so}(N)=\mathfrak{so}_Q(V)=\{x\in\mathfrak{gl}(V)|\la xv,v\ra=0,\ \forall\ v\in V\text{ and }x|_R=0\}$ where  $\la,\ra$ is the bilinear form  associated to $Q$ (namely $\la v,v'\ra=Q(v+v')-Q(v)-Q(v')$ for all $v,v'\in V$) and $R=\{v\in V|\la v,V\ra=0\}$ is the radical of $Q$; recall that $R=0$ unless $p=2$ and $\dim V$ is odd in which case $\dim R=1$ and $Q:R\to\tk$ is injective). 

We have (see \text{\cite{Lu1,LS,Spal2,X1}}) \begin{eqnarray*}
&&\Lambda_{SO(2n+1)}^1=\{(\alpha,\beta)\in\cP_2(n)|\alpha_{i+1}\leq\beta_i\leq\alpha_i+2\},\\
&&\Lambda_{Sp(2n)}^1=\{(\alpha,\beta)\in\cP_2(n)|\alpha_{i+1}-1\leq\beta_i\leq\alpha_i+1\},\\
&&{\Lambda}_{SO(2n)}^1=\{(\alpha,\beta)\in\cP_2(n)|\alpha_{i+1}-2\leq\beta_i\leq\alpha_i\};\\
&&\Lambda_{SO(2n+1)}^2=\Lambda_{Sp(2n)}^2=\{(\alpha,\beta)\in\cP_2(n)|\alpha_{i+1}-2\leq\beta_i\leq\alpha_i+2\},
\\
&&{\Lambda}_{SO(2n)}^2=\{(\alpha,\beta)\in\cP_2(n)|\alpha_{i+1}-4\leq\beta_i\leq\alpha_i\};\\&&
\Lambda_{\mathfrak{so}(2n+1)}^2=\{(\alpha,\beta)\in\cP_2(n)|\beta_i\leq\alpha_i+2\},\ \Lambda_{\mathfrak{sp}(2n)}^2=\cP_2(n),\\
&&{\Lambda}_{\mathfrak{so}(2n)}^2=\{(\alpha,\beta)\in\cP_2(n)|\beta_i\leq\alpha_i\},\end{eqnarray*}
where for $G=SO(2n)$ each pair $(\alpha,\alpha)$ is counted twice in the sets $\Lambda^p$. Note that (see also \cite{Lu5}) $\Lambda_G^1\subset\Lambda_{G}^2\subset\Lambda_{\Lg}^2.$

\subsection{Combinatorial parametrization of $\Omega_G^1$ and the Springer correspondence ($p\neq 2$)}\label{charnot2}Assume $p\neq 2$ and $G=Sp(N)$ (resp. $SO(N)$).
We can identify
\begin{eqnarray*}&&\Omega_{G}^1\text{ with the set }\{\lambda\in\cP(N)|\ m_\lambda(i)\text{ is even if }  i \text{ is odd (resp. even) and }i\neq 0\},
\end{eqnarray*}
where if $G=SO(2n)$  each $\lambda$ with all parts even is counted twice (it corresponds to two (degenerate) classes conjugate under $O(2n)$). For $\rc=\lambda\in\Omega_G^1$, the partition $\lambda$ is given by the sizes of Jordan
blocks of $u-1$ where $u\in\rc$.

Assme $\rc=\lambda=(\lambda_1\geq \lambda_2\geq\cdots)\in{\Omega}^1_G$ and
${\gamma}^1_G(\rc)=(\alpha,\beta)$, $\alpha=(\alpha_1\geq\alpha_2\geq\cdots)$, $\beta=(\beta_1\geq\beta_2\geq\cdots)$.
Recall that $\lambda$ and $(\alpha,\beta)$ are related as follows
\cite{Lu6}. If $G=SO(2n+1)$, then
\begin{equation*}\label{soodd-charnot2}\lambda_{2i-1}=2\alpha_i+1+\delta_i,\quad\lambda_{2i}=2\beta_i-1+\theta_i,
\end{equation*}where
\begin{eqnarray}\label{eqn-o1}
\delta_i=\left\{\begin{array}{lll}1&\text{if
}\beta_i=\alpha_i+2\\-1&\text{if }\alpha_i=\beta_{i-1}\ (i\geq
2)\\0&\text{otherwise}\end{array}\right.,\
\theta_i=\left\{\begin{array}{lll}1&\text{if
}\beta_i=\alpha_{i+1}\\-1&\text{if
}\beta_i=\alpha_{i}+2\\0&\text{otherwise}\end{array}\right.;
\end{eqnarray}
if $G=Sp(2n)$, then
\begin{equation*}\label{sp-charnot2}\lambda_{2i-1}=2\alpha_i+\delta_i,\quad\lambda_{2i}=2\beta_i+\theta_i,
\end{equation*}where
\begin{eqnarray*}
\delta_i=\left\{\begin{array}{lll}1&\text{if
}\beta_i=\alpha_i+1\\-1&\text{if }\alpha_i=\beta_{i-1}+1\ (i\geq
2)\\0&\text{otherwise}\end{array}\right.,\
\theta_i=\left\{\begin{array}{lll}1&\text{if
}\beta_i=\alpha_{i+1}-1\\-1&\text{if
}\beta_i=\alpha_{i}+1\\0&\text{otherwise}\end{array}\right.;
\end{eqnarray*}
if $G=SO(2n)$, then
\begin{equation*}\label{soeven-charnot2}\lambda_{2i-1}=2\alpha_i-1+\delta_i,\quad\lambda_{2i}=2\beta_i+1+\theta_i,\end{equation*}where
\begin{eqnarray*}
\delta_i=\left\{\begin{array}{lll}1&\text{if
}\beta_i=\alpha_i\\-1&\text{if }\alpha_i=\beta_{i-1}+2\ (i\geq
2)\\0&\text{otherwise}\end{array}\right.,\
\theta_i=\left\{\begin{array}{lll}1&\text{if
}\beta_i=\alpha_{i+1}-2\\-1&\text{if
}\beta_i=\alpha_{i}\\0&\text{otherwise}\end{array}\right.
\end{eqnarray*}
(note that ${\gamma}^1_{SO(2n)}(\lambda)=(\alpha,\alpha)$ iff all $\lambda_i$ are even; the two degenerate classes corresponding to such a $\lambda$ are mapped under $\gamma^1_{SO(2n)}$ to the two degenerate elements of $\mathbf{W}^\wedge$ corresponding to $(\alpha,\alpha)$ respectively).

\subsection{Combinatorial parametrization of $\Omega_G^2$ and the Springer correspondence ($p=2$)}\label{sec-1}
Assume $p=2$ in this subsection. Let $G=Sp(2n)=Sp_{\la,\ra}(V)$ (resp. $G=SO(2n)=SO_Q(V)$).  The $Sp(2n)$ (resp. $O(2n)$)-orbit of
$u\in \cU_G$ is characterized by the partition $\lambda\in\cP(2n)$ given by
the sizes of Jordan blocks of $u-1$ and a map
$\varepsilon:\mathbb{N}=\{0,1,2,\ldots\}\rightarrow\{\omega,0,1\}$ satisfying the following conditions (a)-(d) (see \cite[I 2.6]{Spal})\\[5pt]
\indent(a) $\varepsilon(i)=\omega, \text{ if }i\text{ is odd},\text{ or if
}i\geq 1\text{ and }
m_\lambda(i)=0$,\newline
\indent(b) $\varepsilon(i)=1,\text{ if
}i\neq 0\text{ is even and }m_\lambda(i)\text{ is odd},$\\
\indent(c)
$\varepsilon(i)\neq\omega,\text{ if }i\text{ is even and }m_\lambda(i)>0,$\\
\indent(d) $\varepsilon(0)=1\ (\text{resp. }\varepsilon(0)=0).$\\[5pt]
We have  $m_\lambda(i)$ is even for all odd $i$;  $\varepsilon (i)=1$ (for even $i$) iff $\la (u-1)^{\frac{i}{2}}v,(u-1)^{\frac{i}{2}-1}v\ra\neq 0$ (resp. $Q((u-1)^{\frac{i}{2}}v)\neq 0$) for some $v\in V$; $\lambda_1^*$ is even if $G=SO(2n)$.

We identify ${\Omega}_{G}^2$ with the set of all $(\lambda,\varepsilon)$ as above where if $G=SO(2n)$ each $(\lambda,\varepsilon)$ with
$\varepsilon(\lambda_i)=0$ for all $\lambda_i$ is counted twice (it corresponds to two (degenerate) classes conjugate under $O(2n)$). Moreover we identify  $\Omega_{SO(2n+1)}^2$ with $\Omega_{Sp(2n)}^2$ via the natural bijection given by the special isogeny $SO(2n+1)\to Sp(2n)$ (henceforth a class $(\lambda,\varepsilon)\in\Omega_{SO(2n+1)}^2$ is parametrized using $\lambda\in\cP(2n)$ instead of $\lambda\in\cP(2n+1)$).

Assume $\rc=(\lambda,\varepsilon)\in\Omega^2_{G}$, $\lambda=(\lambda_1\geq\lambda_2\geq\cdots)$, and
$\gamma_{G}^2(\rc)=(\alpha,\beta)$, $\alpha=(\alpha_1\geq\alpha_2\geq\cdots)$, $\beta=(\beta_1\geq\beta_2\geq\cdots)$.
Recall
that $(\lambda,\varepsilon)$ and $(\alpha,\beta)$ are related as
follows \cite{LS}. If $G=Sp(2n)$, then
\begin{equation}\label{eqn-4}\lambda_{2i-1}=2\alpha_i+\delta_i,\quad\lambda_{2i}=2\beta_i+\theta_i
,\quad\varepsilon(\lambda_{2i-1})=\varepsilon(\delta_i),
\quad\varepsilon(\lambda_{2i})=\varepsilon(\theta_i),\end{equation}where
\begin{eqnarray}\label{eqn-sp1}
&&\delta_i=\left\{\begin{array}{lllll}2&\text{if
}\beta_i=\alpha_i+2\\1&\text{if }\beta_i=\alpha_i+1\\-2&\text{if
}\alpha_i=\beta_{i-1}+2\ (i\geq 2)\\-1&\text{if
}\alpha_i=\beta_{i-1}+1\ (i\geq
2)\\0&\text{otherwise}\end{array}\right.,\
\theta_i=\left\{\begin{array}{lll}2&\text{if
}\beta_i=\alpha_{i+1}-2\\1&\text{if
}\beta_i=\alpha_{i+1}-1\\-2&\text{if
}\beta_i=\alpha_{i}+2\\-1&\text{if
}\beta_i=\alpha_{i}+1\\0&\text{otherwise}\end{array}\right.,\\
&&\varepsilon(\delta_i)=\left\{\begin{array}{lll}0&\text{if
}\delta_i=\pm2\\ \omega &\text{if
}\delta_i=\pm1\\1&\text{otherwise}\end{array}\right.,\
\varepsilon(\theta_i)=\left\{\begin{array}{lll}0&\text{if }\theta_i=\pm2\\
\omega &\text{if
}\theta_i=\pm1\\1&\text{otherwise}\end{array}\right.;\nonumber
\end{eqnarray}
if $G=SO(2n)$, then\begin{equation*}\label{soeven-char2}\lambda_{2i-1}=2\alpha_i-2+\delta_i\qquad\lambda_{2i}=2\beta_i+2+\theta_i,
\quad\varepsilon(\lambda_{2i-1})=\varepsilon(\delta_i),\quad\varepsilon(\lambda_{2i})=\varepsilon(\theta_i),\end{equation*}where
\begin{eqnarray*}
&&\delta_i=\left\{\begin{array}{lllll}2&\text{if
}\beta_i=\alpha_i\\1&\text{if }\beta_i=\alpha_i-1\\-2&\text{if
}\alpha_i=\beta_{i-1}+4\ (i\geq 2)\\-1&\text{if
}\alpha_i=\beta_{i-1}+3\ (i\geq
2)\\0&\text{otherwise}\end{array}\right.,\
\theta_i=\left\{\begin{array}{lll}2&\text{if
}\beta_i=\alpha_{i+1}-4\\1&\text{if
}\beta_i=\alpha_{i+1}-3\\-2&\text{if
}\beta_i=\alpha_{i}\\-1&\text{if
}\beta_i=\alpha_{i}-1\\0&\text{otherwise}\end{array}\right.,\\
&&\varepsilon(\delta_i)=\left\{\begin{array}{lll}0&\text{if
}\delta_i=\pm2\\ \omega &\text{if
}\delta_i=\pm1\\1&\text{otherwise}\end{array}\right.,\
\varepsilon(\theta_i)=\left\{\begin{array}{lll}0&\text{if }\theta_i=\pm2\\
\omega &\text{if
}\theta_i=\pm1\\1&\text{otherwise}\end{array}\right.
\end{eqnarray*}
(note that $\gamma_{SO(2n)}^2((\lambda,\varepsilon))=(\alpha,\alpha)$ iff
$\varepsilon(\lambda_i)=0$ for all $\lambda_i$; the two degenerate classes corresponding to such a $(\lambda,\varepsilon)$ are mapped under $\gamma^2_{SO(2n)}$ to the two degenerate elements of $\mathbf{W}^\wedge$ corresponding to $(\alpha,\alpha)$ respectively).

\subsection{Combinatorial parametrization of $\Omega_\Lg^2$ and the Springer correspondence ($p=2$)}\label{nilpchar2}Assume $p=2$ and $G=Sp(N)=Sp_{\la,\ra}(V)$ (resp. $SO(N)=SO_Q(V)$).
The $Sp(N)$ (resp. $O(N)$)-orbit
of $x\in\cN_\Lg$ is characterized by the partition
$\lambda=(\lambda_1\geq\lambda_2\geq\cdots)\in\cP(N)$ given by the sizes of Jordan blocks of $x$ and a map $\chi:\{\lambda_i\}_{i\geq 1}\to\mathbb{N}$ satisfying the following conditions (a)-(c) (see \cite{Hes})
\\[5pt]\indent(a) $0\leq\chi(\lambda_i)\leq\frac{\lambda_i}{2}$ (resp.
$[\frac{\lambda_i+1}{2}]\leq\chi(\lambda_i)\leq\lambda_i$),\\
\indent (b) $\chi(\lambda_i)\geq\chi(\lambda_{i+1}),\lambda_i-\chi(\lambda_i)\geq\lambda_{i+1}-\chi(\lambda_{i+1})$,
\\\indent(c) $\chi(\lambda_i)=\frac{\lambda_i}{2}$ (resp. $\chi(\lambda_i)=\lambda_i$), if $m_\lambda(\lambda_i)$ is odd.\\[5pt]
We have $m_\lambda(i)$ is even for all odd $i$ (resp. $\{i\neq 0|m_\lambda(i)\text{ is odd}\}=\{j,j-1\}\cap\mathbb{N}_{+}$ for some $j\in\mathbb{N}_{+}=\{1,2,\ldots\}$); $\chi(a)=\min\{i\in\mathbb{N}|x^{a}v=0\Rightarrow\la x^{2i+1}v,v\ra=0,v\in V\}$ (resp. $\chi(a)=\min\{i\in\mathbb{N}|x^{a}v=0\Rightarrow Q( x^{i}v)=0,v\in V\}$).

 We identify $\Omega_{\Lg}^2$  with the set of all $(\lambda,\chi)$ as above where
if $G=SO(2n)$ each $(\lambda,\chi)$ with $\chi(\lambda_i)=\lambda_i/2$ for all $\lambda_i$ is counted twice (it corresponds to two (degenerate) orbits conjugate under $O(2n)$).

Assume $\rc=(\lambda,\chi)\in\Omega_{\Lg}^2$, $\lambda=(\lambda_1\geq\lambda_2\geq\cdots)$, and
$\gamma_{\Lg}^2(\rc)=(\alpha,\beta)$, $\alpha=(\alpha_1\geq\alpha_2\geq\cdots)$, $\beta=(\beta_1\geq\beta_2\geq\cdots)$.
 Recall
that $(\lambda,\chi)$ and $(\alpha,\beta)$ are related as follows
\cite{Spal2,X1}. If $G=Sp(2n)$, then
\begin{eqnarray*}\label{eqn-3}
&&\lambda_1=\left\{\begin{array}{lll}\alpha_1+\beta_1&\text{ if
}\alpha_1<\beta_1\\
2\alpha_1&\text{ if
}\alpha_1\geq\beta_1\end{array}\right.,\quad\chi(\lambda_1)=\alpha_1,\nonumber\\
&&\lambda_{2i}=\left\{\begin{array}{lll}\alpha_{i+1}+\beta_i&\text{
if }\beta_i<\alpha_{i+1}\\
\alpha_i+\beta_i&\text{ if }\beta_i>\alpha_i\\2\beta_i&\text{ if
}\alpha_{i+1}\leq\beta_i\leq\alpha_i\end{array}\right.,\quad\chi(\lambda_{2i})=\left\{\begin{array}{ll}\alpha_{i}&\text{
if }\beta_i>\alpha_{i}\\
\beta_i&\text{ if }\beta_i\leq\alpha_i\end{array}\right.,\\
&&\lambda_{2i+1}=\left\{\begin{array}{lll}\alpha_{i+1}+\beta_{i}&\text{
if }\alpha_{i+1}>\beta_{i}\\
\alpha_{i+1}+\beta_{i+1}&\text{ if }\alpha_{i+1}<\beta_{i+1}\\2\alpha_{i+1}&\text{ if
}\beta_{i+1}\leq\alpha_{i+1}\leq\beta_{i}\end{array}\right.,\quad\chi(\lambda_{2i+1})=\left\{\begin{array}{ll}\alpha_{i+1}&\text{
if }\alpha_{i+1}\leq\beta_{i}\\
\beta_{i+1}&\text{ if }\alpha_{i+1}>\beta_{i}\end{array}\right.,\ i\geq
1;\nonumber
\end{eqnarray*}
if $G=SO(2n+1)$, let $k\geq 0$ be the largest integer such that $\beta_k>0$, then
\begin{eqnarray*}\label{soodd-nil1}
\lambda_{2i-1}=\left\{\begin{array}{lll}\alpha_{i}+\beta_{i}&\text{
if }i<k+1\\
\alpha_i+1&\text{ if }i=k+1\\ \alpha_i&\text{ if
}i>k+1\end{array}\right.&&\quad\chi(\lambda_{2i-1})=\left\{\begin{array}{ll}\alpha_{i}+1&\text{
if }i\leq k+1\\
\alpha_i&\text{ if }i>k+1\end{array}\right.\\\label{soodd-nil2}\lambda_{2i}=\left\{\begin{array}{lll}\alpha_{i}+\beta_i&\text{
if }i<k+1\\
\alpha_i&\text{ if }i\geq
k+1\end{array}\right.&&\quad\chi(\lambda_{2i})=\left\{\begin{array}{ll}\alpha_{i}+1&\text{
if }i<k+1\\
\alpha_i&\text{ if }i\geq k+1\end{array}\right.,\ i\geq 1;
\end{eqnarray*}
if $G=SO(2n)$, then
\begin{eqnarray*}\label{soeven-nil}
&&\lambda_{2i-1}=\lambda_{2i}=\alpha_i+\beta_i,\quad\chi(\lambda_{2i-1})=\chi(\lambda_{2i})=\alpha_i,\ i\geq 1
\end{eqnarray*}
(note that $\gamma_{\mathfrak{so}({2n})}^2((\lambda,\chi))=(\alpha,\alpha)$ iff
$\chi(\lambda_i)=\lambda_i/2$ for all $i$; the two degenerate classes corresponding to such a $(\lambda,\chi)$ are mapped  under $\gamma^2_{\mathfrak{so}({2n})}$ to the two degenerate elements of $\mathbf{W}^\wedge$ corresponding to $(\alpha,\alpha)$ respectively).

\subsection{Pieces in symplectic groups}\label{lem-sp1}
Assume $G=Sp(2n)$. When $p\neq 2$, each unipotent piece consists of one unipotent class. Let $\rc_i=(\lambda_{\rc_i},\varepsilon_i)\in\Omega_{ G}^2$ (resp.  $\rc_i=(\lambda_{\rc_i},\chi_i)\in\Omega_{\Lg}^2$), $i=1,2$. 
\begin{lem}[\cite{Lu4}]
The classes $\rc_1$ and $\rc_2$ lie in the same unipotent (resp. nilpotent) piece if and only if
$\lambda_{\rc_1}=\lambda_{\rc_2}$.
\end{lem}

\subsection{Pieces in special orthogonal groups}\label{lem-o1}
Assume that $G=SO_Q(V)$ in this subsection. Let  $u\in\cU_G$ (resp.
$x\in \cN_\Lg$). There is  a canonical $Q$-filtration $V_*=(V_{\geq a})_{a\in\mathbb{Z}}$ (where $V_{\geq a+1}\subset V_{\geq
a}\subset V$)
 of $V$ associated to $u$ (resp. $x$) as follows  (see
\cite[2.7(a),  A.4(a)]{Lu4}).
\noindent Let $T=u-1$ (resp. $T=x$). If $p\neq 2$, then

\qquad (a) {\em $V_{\geq a}=\sum_{j\geq\max(0,a)}T^j(\ker T^{2j-a+1})$.}\\[10pt]
If $p=2$, the filtration $V_*=(V_{\geq a})$ is defined by induction on $\dim V$ as follows. If $T=0$ we set $V_{\geq a}=0$ for all $a\geq 1$ and $V_{\geq a}=V$ for all $a\leq0$. Hence $V_*$ is defined when $\dim V\leq 1$. Assume now that $T\neq 0$ and $\dim V\geq 2$. Let $e$ be the smallest integer such that $T^e=0$, $f$ the smallest integer such that $QT^f=0$ and $$m=\max(e-1,2f-2).$$ Then $m\geq 1$. We set
\begin{eqnarray*}&&V_{\geq a}=V\text{ for all } a\leq-m;\ V_{\geq a}=0\text{ for all }a\geq m+1;\\&&V_{\geq -m+1}=\left\{\begin{array}{ll}\{v\in V|T^{e-1}v=0\}&\text{ if }e=2f\\\{v\in V|T^{e-1}v=0, Q(T^{f-1}v)=0\}&\text{ if }e=2f-1\\
\{v\in V|Q(T^{f-1}v)=0\}&\text{ if }e<2f-1\end{array}\right.;\\&&V_{\geq m}=V_{\geq -m+1}^\perp\cap Q^{-1}(0).
\end{eqnarray*}
Let $V'=V_{\geq -m+1}/V_{\geq m}$. Then $Q$ induces a nondegenerate quadratic form $Q'$ on $V'$ and $u$ (resp. $x$) induces a well-defined element $u'\in \cU_{SO_{Q'}(V')}$ (resp. $x'\in\cN_{\mathfrak{so}_{Q'}(V')}$). By induction hypothesis, a canonical $Q'$-filtration $V'_*=(V'_{\geq a})$ of $V'$ is defined for $u'$ (resp. $x'$). For $a\in[-m+1,m]$ we set $V_{\geq a}$ to be the inverse image of $V'_{\geq a}$ under the natural map $V_{\geq -m+1}\to V'$ (note that $V'_{\geq m}=0$ and $V'_{\geq -m+1}=V'$). This completes the definition of $V_*$.

Let $\rc\in\Omega_{G}^p$ (resp.
$\Omega_{\Lg}^2$) and $u\in\rc$ (resp.
$x\in \rc$). Let $V_*=(V_{\geq a})_{a\in\mathbb{Z}}$ be the canonical filtration associated to $u$ (resp. $x$) as above. We
define $$f_a=\dim V_{\geq a}/V_{\geq a+1}.$$
Then 
 $f_a\neq 0$ for finitely many $a\in\mathbb{Z}$ and $f_{-a}=f_a$. 
 The sequence of numbers $(f_a)_{a\in\mathbb{N}}$ ($\mathbb{N}=\{0,1,2,\ldots\}$)  depends only on $\rc$ and not
on the choice of $u\in\rc$ (resp. $x\in\rc$); we denote this sequence  by
$\Upsilon_\rc$. We say two sequences $(f_a)_{a\in\mathbb{N}}=(h_a)_{a\in\mathbb{N}}$ iff $f_a=h_a$ for all $a\in\mathbb{N}$. It follows from  (a) that

\qquad (a$'$) {\em if $\rc=\lambda\in\Omega_G^1$ and $\Upsilon_\rc=(f_a)_{a\in\mathbb{N}}$, then $f_a=\sum_{i\in\mathbb{N}}m_{\lambda}(a+2i+1)\text{ for all }a\in\mathbb{N}$.}

 Let $\rc_1,\rc_2\in\Omega^p$. If $G=SO(2n)$, we assume that $\rc_1$, $\rc_2$ are not  conjugate under $O(2n)$. 
\begin{lem}[\cite{Lu4}]
The classes $\rc_1$ and $\rc_2$ lie in the same piece if and only if
$\Upsilon_{\rc_1}=\Upsilon_{\rc_2}$.
\end{lem}

\subsection{Closure relations on unipotent classes and nilpotent orbits}\label{sec-ord}Let $\rc,\rc'\in\Omega^{p}_{G}$ (resp. $\Omega^{2}_{\Lg}$). We say that
$\rc\leq\rc'$ if $\rc$  is contained in the closure of $\rc'$ in $G$ (resp. $\Lg$); and that $\rc<\rc'$ if $\rc\leq\rc'$ and $\rc\neq\rc'$. In the following if $G=SO(2n)$, we assume that $\rc$ and $\rc'$
are not conjugate under $O(2n)$ (otherwise they are incomparable with respect to the partial order $\leq$).

Assume $\rc=\lambda,\rc'=\lambda'\in\Omega_G^1$. We have $\rc\leq\rc'$ if and only if $\lambda\leq\lambda'$ (see \cite[II
8.2]{Spal}).

Assume $\rc=(\lambda,\varepsilon),\rc'=(\mu,\phi)\in\Omega_{G}^2$. We order the set $\{\omega,0,1\}$ by $\omega<0<1$.
Then $\rc\leq \rc'$ if and only if
$(\lambda,\varepsilon)\leq(\mu,\phi)$ (see \cite[II
8.2]{Spal}), namely,  the following conditions (a)-(c) hold\\[5pt]
\indent (a) $\lambda\leq\mu$,\\[5pt]
\indent (b) $\sum_{j\in[1,i]}\lambda_j^*-\max(\varepsilon(i),0)\geq\sum_{j\in[1,i]}\mu_j^*-\max(\phi(i),0)$, for all $i\geq1$,\\[5pt]
\indent (c) if $\sum_{j\in[1,i]}\lambda_j^*=\sum_{j\in[1,i]}\mu_j^*$ and
$\lambda_{i+1}^*-\mu_{i+1}^*$ is odd then $\phi(i)\neq 0$, for all $i\geq1$.

\section{Reformulation of closure relations on $\Omega^p$}
\subsection{}\label{prop-order}Let $\rc,\rc'\in\Omega^p$
(if $G=SO(2n)$,  we assume that $\rc$, $\rc'$ are not conjugate under $O(2n)$).
\begin{prop}\label{prop-ord}
We have $\rc\leq\rc'$ if and only if
$\gamma^p(\rc)\leq\gamma^p(\rc')$.
\end{prop}
\noindent If $\rc,\rc'\in\Omega_{\Lg}^2$, the proposition
is a result of Spaltenstein \cite{Spal2}.  The proofs for $\rc,\rc'\in\Omega_{G}^p$ when $p\neq 2$ and $p=2$ are given in subsections \ref{sec-order1} and  \ref{sec-order2} respectively.

\subsection{}\label{sec-order1}
Assume $\rc=\lambda,\rc'=\lambda'$ in $\Omega^1_G$,
$\gamma^1_G(\rc)=(\alpha,\beta)$ and
$\gamma^1_G(\rc)=(\alpha',\beta')$. We show that

 \qquad (a) {\em $\lambda\leq\lambda'$ iff
$(\alpha,\beta)\leq(\alpha',\beta')$.}\\[10pt]
We prove (a) for $G=SO(2n+1)$. The proofs for $Sp(2n)$ and $SO(2n)$ are entirely
similar and omitted. For $(\alpha,\beta)\in\Lambda^1_{SO(2n+1)}$, let $A_i,B_i$ be as in (\ref{eqn-ab}) and let
$\Delta_i=\sum_{j\in[1,i]}(\delta_j+\theta_j)$,
$\Theta_i=\Delta_{i-1}+\delta_i$, where
$\delta_j,\theta_j$ are as in (\ref{eqn-o1}).  One can easily verify that

$(*)\ \begin{array}{l}\label{eqn-1}
\Delta_i=\left\{\begin{array}{ll}1&\text{if
}\beta_i=\alpha_{i+1}\\0&\text{otherwise}\end{array}\right.,\
\Theta_i=\left\{\begin{array}{ll}1&\text{if
}\beta_i=\alpha_{i}+2\\0&\text{otherwise}\end{array}\right.
\end{array}$.\\[5pt]
We have $\sum_{j=1}^{2i}\lambda_j=2A_i+\Delta_i$ and
$\sum_{j=1}^{2i-1}\lambda_j=2B_i+\Theta_i+1$. Assume $\lambda\leq\lambda'$. It follows from
\ref{sec-partition} (a) and $(*)$
that $A_i\leq A_i'$ and $B_i\leq B_i'$. Hence
$(\alpha,\beta)\leq(\alpha',\beta')$ (see \ref{sec-partition}). Conversely assume $(\alpha,\beta)\leq(\alpha',\beta')$. Then $A_i\leq A_i'$ and
$B_i\leq B_i'$ for all $i$. We show that $A_i=A_i'$ implies
$\Delta_i\leq\Delta_i'$. Assume otherwise,
$A_i=A_i',\Delta_i=1,\Delta_i'=0$. Then $\beta_i=\alpha_{i+1}$ and
$\beta_i'>\alpha_{i+1}'$. Since $B_{i}=A_i-\beta_i\leq
B_{i}'=A_i'-\beta_i'$ and $B_{i+1}=A_i+\alpha_{i+1}\leq
B_{i+1}'=A_i'+\alpha_{i+1}'$, we have $\beta_i\geq\beta_i'$ and
$\alpha_{i+1}\leq\alpha_{i+1}'$ which is a contradiction. Similarly $B_i=B_i'$ implies $\Theta_i\leq\Theta_i'$. Hence $\lambda\leq\lambda'$.

\begin{rem}
When $G$ is $Sp(2n)$ or $SO(2n+1)$, (a) is also obtained in \cite[Proposition 2.4 and Proposition 2.11]{AHS}. I am grateful to the referee for pointing this out.
\end{rem}

\subsection{}\label{sec-order2}Assume $\rc=(\lambda,\varepsilon),\rc'=(\mu,\phi)\in\Omega_{G}^2$,
$\gamma_{G}^2(\rc)=(\alpha,\beta)$ and $\gamma_{G}^2(\rc')=(\alpha',\beta')$. We show that

\qquad (a) {\em  $(\lambda,\varepsilon)\leq(\mu,\phi)$ iff
$(\alpha,\beta)\leq(\alpha',\beta')$.}\\[10pt]
We prove (a) for $G=Sp(2n)$ (and thus for $G=SO(2n+1)$). The proof for $SO(2n)$ is entirely
similar and omitted. Since
$\sum_{j>i}\lambda_j^*=\sum_{j\in[1,\lambda_{i+1}^*]}(\lambda_j-i)$
and, for $i$ large enough, $\sum_{j\in[1,i]}\lambda_j^*=\sum_{j\in[1,i]}\mu_j^*$, we have

\qquad (b) {\em $\sum_{j\in[1,i]}\lambda_j^*=\sum_{j\in[1,i]}\mu_j^*$ iff $\sum_{j\in[1,\lambda_{i+1}^*]}(\lambda_j-i)=\sum_{j\in[1,\mu_{i+1}^*]}(\mu_j-i)$.}\\[10pt]
We show that

\qquad(c) {\em  if $\lambda\leq\mu$ and $\sum_{j\in[1,k]}\lambda_j^*=\sum_{j\in[1,k]}\mu_j^*$, then   $\sum_{j\in[1,\lambda_{k+1}^*]}\lambda_j=\sum_{j\in[1,\lambda_{k+1}^*]}\mu_j$, }

\qquad(d) {\em if $\lambda\leq\mu$ and $\sum_{j\in[1,m]}\lambda_j=\sum_{j\in[1,m]}\mu_j$, then $\sum_{j\in[1,\mu_m]}\lambda_j^*=\sum_{j\in[1,\mu_m]}\mu_j^*$.}\\[10pt]
By \ref{sec-partition} (a$'$), the assumptions in (c) imply that $\lambda_k^*\leq\mu_k^*$, $\lambda_{k+1}^*\geq\mu_{k+1}^*$. It follows that $\mu_j=k$ for $j\in[\mu_{k+1}^*+1,\lambda_{k+1}^*]$ and thus $\sum_{j\in[1,\mu_{k+1}^*]}(\mu_j-k)=\sum_{j\in[1,\lambda_{k+1}^*]}(\mu_j-k)$. Now (c) follows from (b). By \ref{sec-partition} (a), the assumptions in (d) imply that $\lambda_{m+1}\leq\mu_{m+1}$ and $\lambda_m\geq\mu_m$. Let $k=\mu_m$.  Then $\sum_{j\in[1,\lambda_k^*]}(\lambda_j-k)=\sum_{j\in[1,m]}(\lambda_j-k),\ \sum_{j\in[1,\mu_k^*]}(\mu_j-k)=\sum_{j\in[1,m]}(\mu_j-k)$ (since $\lambda_i=k$, $i\in[m+1,\lambda_k^*]$; $\mu_i=k$, $i\in[m+1,\mu_k^*]$). Now (d) follows from (b).

For $(\alpha,\beta)\in\Lambda_{Sp(2n)}^2$, let $A_i,B_i$ be as in (\ref{eqn-ab}) and let
$\Delta_i=\sum_{j\in[1,i]}(\delta_j+\theta_j)$,
$\Theta_i=\Delta_{i-1}+\delta_i$, where $\delta_j$
and $\theta_j$ are as in (\ref{eqn-sp1}). One can easily verify  that
we have

$\quad(*)\quad
\Delta_i=\left\{\begin{array}{lll}2&\text{if
}\beta_i=\alpha_{i+1}-2\\1&\text{if
}\beta_i=\alpha_{i+1}-1\\0&\text{otherwise}\end{array}\right.,\
\Theta_i=\left\{\begin{array}{lll}2&\text{if
}\beta_i=\alpha_{i}+2\\1&\text{if
}\beta_i=\alpha_{i}+1\\0&\text{otherwise}\end{array}\right.$.\\[5pt]
Using (\ref{eqn-4}) one can easily check that $\lambda_{2i-1}=\lambda_{2i}$ iff $\beta_i\geq\alpha_i$, or (if $i\geq 2$) $\beta_{i}=\beta_{i-1}$, $\alpha_i=\alpha_{i+1}$ and $\beta_i\leq\alpha_{i+1}-1$; and $\lambda_{2i}=\lambda_{2i+1}$ iff $\beta_i\leq\alpha_{i+1}$, or $\beta_{i}=\beta_{i+1}$, $\alpha_i=\alpha_{i+1}$ and $\beta_i\geq\alpha_{i}+1$. It then follows that

\qquad(e) {\em  if $\beta_i=\alpha_i+2$, then $\lambda_{\lambda_{2i}+1}^*$ is even,}

\qquad(f) {\em  if $\alpha_i=\beta_{i-1}+2$ ($i\geq 2$), then $\lambda_{\lambda_{2i-1}+1}^*$ is odd.}\\[10pt]
We have
$\sum_{j\in[1,2i]}\lambda_j=2A_i+\Delta_i$ and
$\sum_{j\in[1,2i-1]}\lambda_j=2B_i+\Theta_i$.

Assume $(\lambda,\varepsilon)\leq(\mu,\phi)$.  It follows from
 $\lambda\leq\mu$ and $(*)$ that $A_i\leq A_i'$ except if
$\Delta_i=0,\Delta_{i}'=2$ and
$\sum_{j\in[1,2i]}\lambda_j=\sum_{j\in[1,2i]}\mu_j$. In the latter case,
we have $\beta_i\geq\alpha_{i+1}$, $\beta_{i}'=\alpha_{i+1}'-2$, $\lambda_{2i+1}\leq \mu_{2i+1}$ and $\lambda_{2i}\geq \mu_{2i}$.
Then $\mu_{2i}=\mu_{2i+1}$ and $\phi(\mu_{2i})=0$ (we use (\ref{eqn-4})). Let
$k=\mu_{2i}$.  By (d), we have
$\sum_{j\in[1,k]}\lambda_j^*=\sum_{j\in[1,k]}\mu_j^*$. By (f), $\mu_{k+1}^*$ is odd. If $\lambda_{2i}>k$, then $\lambda_{k+1}^*=2i$ is even, which contradicts \ref{sec-ord} (c). Hence $\lambda_{2i}=k$ and thus $\varepsilon(k)=0$ (we use  \ref{sec-ord} (b) and $k$ even). It follows that $\beta_i=\alpha_i+2$ (since $\beta_i\geq \alpha_{i+1}$) and thus $\lambda_{k+1}^*$ is even by (e), which again contradicts \ref{sec-ord} (c). Hence $A_i\leq A_i'$. Similarly we have $B_i\leq B_i'$. Hence
$(\alpha,\beta)\leq(\alpha',\beta')$.

Assume $(\alpha,\beta)\leq(\alpha',\beta')$. We show that if
$A_i=A_i'$ then $\Delta_i\leq\Delta_i'$. Assume otherwise,
$A_i=A_i'$, $\Delta_i=1$ (resp. $2$), $\Delta_i'=0$ (resp. $1$ or
$0$). Then as in the proof of \ref{sec-order1} (a), we have $\beta_i\geq\beta_i'$ and
$\alpha_{i+1}\leq\alpha_{i+1}'$, which contradicts to $\Delta_i\leq\Delta_i'$ (we use $(*)$). Similarly one can show if $B_i=B_i'$
then $\Theta_i\leq\Theta_i'$. It follows that $\lambda\leq\mu$.

We verify \ref{sec-ord} (b). Assume $\varepsilon(k)=1$,
$\phi(k)\leq 0$, and
$\sum_{j\in[1,k]}\lambda_j^*=\sum_{j\in[1,k]}\mu_j^*$. Let $m=\lambda_{k+1}^*$. Then $\lambda_{m+1}=\mu_{m+1}=k$ (since $\mu_{k+1}^*\leq\lambda_{k+1}^*<\lambda_k^*\leq\mu_k^*$). By (c), we
have $\sum_{j\in[1,m]}\lambda_j=\sum_{j\in[1,m]}\mu_j$. Suppose $m=2i$. Note $\varepsilon(\lambda_{m+1})=1$
implies that $\delta_{i+1}=0$, $\Delta_i=0$ and thus
$\lambda_{m+1}=2\alpha_{i+1}$.  Since $A_i\leq A_i'$ and $2A_i+\Delta_i=2A_i'+\Delta_i'$, we have
$A_i=A_i'$ and $\Delta_i'=0$. Together with $\phi(\mu_{m+1})\leq 0$,
this implies that $\mu_{m+1}>2\alpha_{i+1}'$. Hence
$\alpha_{i+1}>\alpha_{i+1}'$, which contradicts $B_{i+1}\leq
B_{i+1}'$. Suppose $m=2i-1$.
Note $\varepsilon(\lambda_{m+1})=1$ implies that $\theta_{i}=0$,
$\Theta_i=0$ and thus $\lambda_{m+1}=2\beta_{i}$. Since $B_i\leq
B_i'$ and $2B_i+\Theta_i=2B_i'+\Theta_i'$, we have $B_i=B_i'$ and $\Theta_i'=0$. Together with
$\phi(\mu_{m+1})\leq 0$ this implies that $\mu_{m+1}>2\beta_{i}'$.
Hence $\beta_{i}>\beta_{i}'$, which contradicts $A_{i}\leq A_{i}'$.

It remains to verify \ref{sec-ord} (c). Assume $\sum_{j\in[1,k]}\lambda_j^*=\sum_{j\in[1,k]}\mu_j^*$,
$\lambda_{k+1}^*-\mu_{k+1}^*$ is odd,
and $\phi(k)=0$. Let $\lambda_{k+1}^*=m$. Then $\mu_{m}=k$ (since $\mu_{k+1}^*<\lambda_{k+1}^*\leq \lambda_k^*\leq\mu_k^*$). By (c),
$\sum_{j\in[1,m]}\lambda_j=\sum_{j\in[1,m]}\mu_j$. Suppose $m=2i$. Note that $\varepsilon(\mu_{m})=0$ implies that $\theta_i'=2,\Delta_i'=2$
$(\beta_i'=\alpha_{i+1}'-2)$ or
$\theta_i'=-2,\Delta_i'=0\ (\beta_i'=\alpha_i'+2)$. If
$\beta_i'=\alpha_{i+1}'-2$, then $\Delta_i'=2$ and
$2A_i+\Delta_i=2A_i'+\Delta_i'$ imply that
$A_i=A_i'$, $\Delta_i=2$ and thus $\beta_i=\alpha_{i+1}-2$,
$\theta_i=2$. Since $\lambda_m=2\beta_i+\theta_i>\mu_m=2\beta_i'+\theta_i'$, we have $\beta_i>\beta_i'$ and thus $\alpha_{i+1}>\alpha_{i+1}'$, which contradicts $B_{i+1}\leq B_{i+1}'$. If
$\beta_i'=\alpha_{i}'+2$, then $\mu_{k+1}^*$ is even (see (e)), which
contradicts the fact that $\lambda_{k+1}^*-\mu_{k+1}^*$ is odd.
Suppose $m=2i-1$. Note that
$\phi(\mu_m)=0$ implies that
$\delta_i'=2$, $\Theta_i'=2$ ($\beta_i'=\alpha_i'+2$) or
$\delta_i'=-2$, $\Theta_i'=0$ ($\alpha_i'=\beta_{i-1}'+2$). If
$\beta_i'=\alpha_i'+2$, then $\Theta_i'=2$ and $2B_i+\Theta_i=2B_i'+\Theta_i'$ imply $\Theta_i=2$, $B_i=B_i'$ and thus
$\beta_i=\alpha_i+2$, $\delta_i=2$. Since $\lambda_{m}=2\alpha_i+\delta_i>\mu_{m}=2\alpha_i'+\delta_i'$, we have $\alpha_i>\alpha_i'$ and thus $\beta_i>\beta_i'$, which
contradicts $A_i\leq A_i'$. If
$\alpha_i'=\beta_{i-1}'+2$, then $\mu_{k+1}^*$ is odd, which
contradicts the fact that $\lambda_{k+1}^*-\mu_{k+1}^*$ is odd. This completes the proof of (a).

\section{Combinatorial definition of unipotent and nilpotent pieces}\label{sec-comdef}

\subsection{}\label{ssec-comb1}  Let $\tilde{\rc}\in\Omega_G^1$ and let $\rc\in\Omega_G^p$ (resp. $\rc\in\Omega_\Lg^2$) be such that $\gamma_G^p(\rc)=\gamma^1_G(\tilde{\rc})$ (resp. $\gamma_\Lg^2(\rc)=\gamma^1_G(\tilde{\rc})$). Define
$\Sigma_{\tilde{\rc}}^{p,G}$ (resp. $\Sigma_{\tilde{\rc}}^{2,\Lg}$)
to be the set of all classes $\rc'\in\Omega_G^p$ (resp. $\rc'\in\Omega_\Lg^2$) such that $\rc'\leq \rc$
and
$\rc'\nleq\rc''$ for any
$\rc''<\rc$  with
$\gamma_G^p(\rc'')\in\Lambda_G^1$ (resp. $\gamma_\Lg^2(\rc'')\in\Lambda_G^1$). We show that (see \ref{ssec-fiber1} (a))

(a) {\em$\{\Sigma_{\tilde{\rc}}^{p,G}\}_{\tilde{\rc}\in\Omega_G^1}$ (resp. $\{\Sigma_{\tilde{\rc}}^{2,\Lg}\}_{\tilde{\rc}\in\Omega_G^1}$)
form a partition of $\cU_G$ (resp. $\cN_\Lg$) }\\[8pt]
and that (see \ref{ssec-map})

 (b) {\em each set $\Sigma_{\tilde{\rc}}^{p,G}$ (resp. $\Sigma_{\tilde{\rc}}^{2,\Lg}$), ${\tilde{\rc}\in\Omega_G^1}$, is a unipotent (resp. nilpotent) piece defined in \cite{Lu4}.}\\[10pt]
In view of the properties of unipotent (resp. nilpotent) pieces proved by Lusztig (resp. and the author) \cite{Lu4}, (b) implies  that each set $\Sigma_{\tilde{\rc}}^{p,G}$ (resp. $\Sigma_{\tilde{\rc}}^{2,\Lg}$) is smooth and that the number of $\tF_{p^s}$-rational points in such a set (called an MS-piece by Lusztig in unipotent case) is a polynomial in $p^s$ with integer coefficients independent of $p,s$ (in unipotent case this is the statement 6.8 (a) of \cite{Lu7} for classical groups).  The definition of unipotent pieces using closure relations  is first considered by Spaltenstein and (a) for $\cU_G$ is shown in \cite{Spal}. For completeness, we include here a different proof of (a) which applies for both $\cU_G$ and $\cN_\Lg$.

We define maps
$$ \Phi_G: \Lambda_\Lg^2\rightarrow\Lambda_G^1,\
(\alpha,\beta)\mapsto(\tilde{\alpha},\tilde{\beta})$$
as follows. We can restrict $\Phi_G$ to the sets $\Lambda_G^1$ and $\Lambda_G^2$ using the inclusion $\Lambda_G^1\subset\Lambda_G^2\subset\Lambda_\Lg^2$.\\[10pt]
If $G=SO(2n+1)$, define $\tilde{\alpha}_1=\alpha_1$,
\begin{eqnarray*}\label{phi-soodd1}
&&\tilde{\alpha}_i=\left\{\begin{array}{ll}[\frac{\alpha_i+\beta_{i-1}}{2}]&\text{
if }\alpha_i>\beta_{i-1}\\ \alpha_i&\text{ if
}\alpha_i\leq\beta_{i-1}\end{array}\right.i\geq 2,\quad
\tilde{\beta}_i=\left\{\begin{array}{ll}[\frac{\alpha_{i+1}+\beta_{i}+1}{2}]&\text{
if }\beta_i<\alpha_{i+1}\\ \beta_i&\text{ if
}\beta_i\geq\alpha_{i+1}\end{array}\right.i\geq1.
\end{eqnarray*}
If $G=Sp(2n)$, define $\tilde{\alpha}_1=\left\{\begin{array}{ll}[\frac{\alpha_1+\beta_1}{2}]&\text{
if }\beta_1>\alpha_1+1\\ \alpha_1&\text{ if
 }\beta_1\leq\alpha_1+1\end{array}\right.$ and
\begin{eqnarray*}\label{phi-sp1}
\tilde{\alpha}_i=\left\{\begin{array}{ll}[\frac{\alpha_i+\beta_{i}}{2}]\quad \text{if }\alpha_{i}<\beta_i-1\\
{}[\frac{\alpha_i+\beta_{i-1}+1}{2}]\quad\text{if }\alpha_i>\beta_{i-1}+1\\
\alpha_i\quad \text{if }\beta_i-1\leq\alpha_{i}\leq
\beta_{i-1}-1\end{array}\right. i\geq 2,\ \tilde{\beta}_i=\left\{\begin{array}{lll}[\frac{\alpha_{i}+\beta_{i}+1}{2}]&\text{if }\beta_i>\alpha_{i}+1\\{}[\frac{\alpha_{i+1}+\beta_{i}}{2}]&\text{if }\beta_i<\alpha_{i+1}-1\\
\beta_i&\text{if
}\alpha_{i+1}-1\leq\beta_i\leq\alpha_{i}+1\end{array}\right.i\geq1.
\end{eqnarray*}
If $G=SO(2n)$, define $\tilde{\alpha}_1=\alpha_1$ and
\begin{eqnarray*}\label{phi-soeven1}
\tilde{\alpha}_i=\left\{\begin{array}{ll}[\frac{\alpha_i+\beta_{i-1}+2}{2}]&\text{
if }\alpha_i>\beta_{i-1}+2\\ \alpha_i&\text{ if
}\alpha_i\leq\beta_{i-1}+2\end{array}\right.\ i\geq 2,
\ \tilde{\beta}_i=\left\{\begin{array}{ll}[\frac{\alpha_{i+1}+\beta_{i}-1}{2}]&\text{
if }\beta_i<\alpha_{i+1}-2\\ \beta_i&\text{ if
}\beta_i\geq\alpha_{i+1}-2\end{array}\right.\ i\geq1,
\end{eqnarray*}
(note that $\Phi_G((\alpha,\beta))=(\tilde{\alpha},\tilde{\alpha})$ iff $(\alpha,\beta)=(\tilde{\alpha},\tilde{\alpha})$); we define $\Phi_G$ to be the identity map on the set of  degenerate elements of $\mathbf{W}^\wedge$.

It is easy to verify that in each case we get a well-defined element
$(\tilde{\alpha},\tilde{\beta})\in \Lambda^1_G$.

\subsection{}\label{ssec-fiber1}
In this subsection we show that for each $\tilde{\rc}\in\Omega_G^1$,

\qquad(a) {\em$\gamma_G^p({\Sigma}_{\tilde{\rc}}^{p,G})=(\Phi_G|_{\Lambda_G^p})^{-1}(\gamma_G^1(\tilde{\rc}))$ (resp. $\gamma_\Lg^2({\Sigma}_{\tilde{\rc}}^{2,\Lg})=\Phi_G^{-1}(\gamma_G^1(\tilde{\rc}))$).}\\[10pt]
Then  \ref{ssec-comb1} (a) follows from (a).  In view of Proposition \ref{prop-order}, (a) follows from the definition of ${\Sigma}_{\tilde{\rc}}^{p,G}$ (resp. ${\Sigma}_{\tilde{\rc}}^{2,\Lg}$) and the following (b) and (c)

 \qquad(b) {\em $\Phi_G|_{\Lambda_G^1}=Id$ and $(\alpha,\beta)\leq\Phi_G((\alpha,\beta))$,}

\qquad(c) {\em For any $(\tilde{\alpha}',\tilde{\beta}')\in\Lambda_G^1$ such that $(\alpha,\beta)\leq(\tilde{\alpha}',\tilde{\beta}')$, we have $\Phi_G((\alpha,\beta))\leq(\tilde{\alpha}',\tilde{\beta}')$.}
\begin{rem}
When $G=Sp(2n)$ (resp. $SO(2n+1)$), (b) and (c) are also obtained in Corollary 2.6 (resp. Corollary 2.17) of \cite{AHS}. In fact, the image $\Phi_G((\alpha,\beta))$ of $(\alpha,\beta)\in\Lambda^2_\Lg$ under our map $\Phi_G$ is the same as $(\alpha,\beta)^C$ (resp. $(\alpha,\beta)^B$) in the notation of \cite{AHS}. I am grateful to the referee for pointing this out. \end{rem}
The first assertion in (b) follows from the definition of $\Phi_G$. Suppose $\Phi_G((\alpha,\beta))=(\tilde{\alpha},\tilde{\beta})$. Let $A_i$,
$B_i$,  $\tilde{A}_i,\tilde{B}_i$, $\tilde{A}_i',\tilde{B}_i'$ be defined for $(\alpha,\beta)$, $(\tilde{\alpha},\tilde{\beta})$, $(\tilde{\alpha}',\tilde{\beta}')$ respectively as in (\ref{eqn-ab}). We prove the second assertion in (b) and (c) in various cases.

(i) Assume $G=SO(2n+1)$. Note that we have
$\tilde{\beta}_i+\tilde{\alpha}_{i+1}={\beta}_i+{\alpha}_{i+1}$,
$B_1=\tilde{B}_1$ and thus $B_i=\tilde{B}_i$. Moreover, $A_i\leq\tilde{A}_i$, and $A_i<\tilde{A}_i$ if and only if $\beta_i<\alpha_{i+1}$. Hence
$(\alpha,\beta)\leq(\tilde{\alpha},\tilde{\beta})$.

Assume there exists
$(\tilde{\alpha}',\tilde{\beta}')\in\Lambda_{G}^1$ such
that
$(\alpha,\beta)\leq(\tilde{\alpha}',\tilde{\beta}')$
and $(\tilde{\alpha},\tilde{\beta})\nleq(\tilde{\alpha}',\tilde{\beta}')$. Since
$\tilde{B}_j=B_j\leq\tilde{B}_j'$ for all $j$, there exists an $i$
such that $\tilde{A}_i'<\tilde{A}_i$. It follows that $\beta_i<\alpha_{i+1}$ (since $A_i<\tilde{A}_i$) and thus $\tilde{\beta}_i\leq\tilde{\alpha}_{i+1}+1$ by the definition of $\Phi_G$. On the other hand,
$\tilde{\beta}_i>\tilde{\beta}_i'\geq\tilde{\alpha}_{i+1}'>\tilde{\alpha}_{i+1}$ (we use $\tilde{B}_j=B_j\leq\tilde{B}_j'$, $j=i,i+1$, and the fact that $(\tilde{\alpha}',\tilde{\beta}')\in\Lambda_{G}^1$), which is a contradiction.

(ii) Assume $G=Sp(2n)$. We show by induction on $i$ that

\quad (d) {\em if
$\beta_i>\alpha_{i}+1$, then $\tilde{B}_i>B_i$, $\tilde{A}_i=A_i$;
if $\alpha_{i+1}-1\leq\beta_i\leq\alpha_{i}+1$, then
$\tilde{A}_i=A_i$, $\tilde{B}_i=B_i$; if $\beta_i<\alpha_{i+1}-1$,
then $\tilde{B}_i=B_i$, $\tilde{A}_i>A_i$.}\\[10pt]
It then follows that
$(\alpha,\beta)\leq(\tilde{\alpha},\tilde{\beta})$.
It is easy to verify that
(d) holds when $i=1$. We have the following subcases:

(ii-1) $\beta_{i+1}>\alpha_{i+1}+1$. Then $\tilde{\alpha}_{i+1}>\alpha_{i+1}$ and $\tilde{\alpha}_{i+1}+\tilde{\beta}_{i+1}={\alpha}_{i+1}+{\beta}_{i+1}$. Since
$\beta_{i}>\alpha_{i+1}+1$, by induction hypothesis, $\tilde{A}_i=A_i$. It follows that
$\tilde{B}_{i+1}>\tilde{B}_{i+1}$ and $\tilde{A}_{i+1}=A_{i+1}$.

(ii-2)
$\alpha_{i+2}-1\leq\beta_{i+1}\leq\alpha_{i+1}+1$. If
$\beta_i\geq\alpha_{i+1}-1$, then $\tilde{A}_i={A}_i$ (by induction hypothesis) and
$\tilde{\alpha}_{i+1}={\alpha}_{i+1}$;
if $\beta_i<\alpha_{i+1}-1$, then $\tilde{B}_i={B}_i$ (by induction hypothesis) and
$\tilde{\beta}_i+\tilde{\alpha}_{i+1}={\beta}_i+{\alpha}_{i+1}$. It follows that $\tilde{B}_{i+1}=B_{i+1}$.
Since $\tilde{\beta}_{i+1}=\beta_{i+1}$, we have
$\tilde{A}_{i+1}=A_{i+1}$.

(ii-3) $\beta_{i+1}<\alpha_{i+2}-1$. If $\beta_i\geq\alpha_{i+1}-1$,
then $\tilde{A}_i=A_i$ and $\tilde{\alpha}_{i+1}=\alpha_{i+1}$; if
$\beta_i<\alpha_{i+1}-1$, then $\tilde{B}_i=B_i$ and
$\tilde{\beta}_i+\tilde{\alpha}_{i+1}={\beta}_i+\alpha_{i+1}$. It
follows that $\tilde{B}_{i+1}=B_{i+1}$. Since
$\tilde{\beta}_{i+1}>\beta_{i+1}$, we have $\tilde{A}_{i+1}>A_{i+1}$. (d) is proved.

Assume there exists
$(\tilde{\alpha}',\tilde{\beta}')\in\Lambda_{G}^1$ such
that
$(\alpha,\beta)\leq(\tilde{\alpha}',\tilde{\beta}')$
and $(\tilde{\alpha},\tilde{\beta})\nleq(\tilde{\alpha}',\tilde{\beta}')$. Suppose
that there exists an $i$ such that
$\tilde{A}_i'<\tilde{A}_i$. Then it follows from (d) that $\beta_i<\alpha_{i+1}-1$ (since $A_i<\tilde{A}_i$) and thus $\tilde{\beta}_i\leq\tilde{\alpha}_{i+1}$; $\tilde{B}_i=B_i$, $\tilde{B}_{i+1}=B_{i+1}$ and thus
$\tilde{\beta}_i\geq\tilde{\beta}_i'+1\geq\tilde{\alpha}_{i+1}'\geq\tilde{\alpha}_{i+1}+1$, which is a contradiction.
Then there exists an $i$ such that
$\tilde{B}_i'<\tilde{B}_i$. It follows from (d) that  $\beta_i>\alpha_{i}+1$ and thus $\tilde{\beta}_i\geq\tilde{\alpha}_{i}$; $\tilde{A}_i=A_i\leq\tilde{A}_i'$,
$\tilde{A}_{i-1}=A_{i-1}\leq\tilde{A}_{i-1}'$, and thus
$\tilde{\beta}_i\leq\tilde{\beta}_i'-1\leq\tilde{\alpha}_{i}'\leq\tilde{\alpha}_{i}-1$, which is again a contradiction.

(iii) Assume $G=SO(2n)$. We have
$\tilde{\beta}_i+\tilde{\alpha}_{i+1}={\beta}_i+{\alpha}_{i+1}$,
$B_1=\tilde{B}_1$ and thus $B_i=\tilde{B}_i$. Moreover, $A_i\leq\tilde{A}_i$, and $A_i<\tilde{A}_i$ if and only if  $\beta_i<\alpha_{i+1}-2$. Hence
$(\alpha,\beta)\leq(\tilde{\alpha},\tilde{\beta})$.

Assume there exists
$(\tilde{\alpha}',\tilde{\beta}')\in\Lambda_{G}^1$ such
that
$(\alpha,\beta)\leq(\tilde{\alpha}',\tilde{\beta}')$
and $(\tilde{\alpha},\tilde{\beta})\nleq(\tilde{\alpha}',\tilde{\beta}')$. Then
$\tilde{B}_j=B_j\leq\tilde{B}_j'$ for all $j$ and there exists an $i$
such that $A_i\leq \tilde{A}_i'<\tilde{A}_i$. It follows that $\tilde{\beta}_i<\tilde{\alpha}_{i+1}$, and
$\tilde{\beta}_i>\tilde{\beta}_i'\geq \tilde{\alpha}_{i+1}'-2>\tilde{\alpha}_{i+1}-2$, which is
a contradiction. This completes the proof of (b) and (c).

\section{Explicit description of pieces}\label{sec-fiber}

\subsection{}\label{ssec-map} We define  maps $$\Psi_G^p:\Omega_{G}^p\rightarrow\Omega_{G}^1 \text{ and  }\Psi_\Lg^2:\Omega_{\Lg}^2\rightarrow \Omega_{G}^1$$ as follows such that their fibers are pieces (see Proposition \ref{prop-psi}). When $p\neq 2$ let $\Psi_G^p$ be the natural identification map between $\Omega_G^p$ and $\Omega_G^1$.

Assume that $G=Sp(2n)$. Define $\Psi_G^2((\lambda,\varepsilon))=\lambda\text{ and } \Psi_\Lg^2((\lambda,\chi))=\lambda$.

Assume that $G=SO(N)$ and $\rc=(\lambda,\varepsilon)\in\Omega_G^2$ (note that if $N$ is odd we use the identification $\Omega_G^2=\Omega_{Sp(N-1)}^2$ and thus $\lambda\in\cP(N-1)$). Define  $\Psi_G^2(\rc)=\tilde{\lambda}$ as follows.
If $\varepsilon(\lambda_{2i-1})=1$ and  (when $i\geq 2$) $ \lambda_{2i-1}<\lambda_{2i-2}$, then $\tilde{\lambda}_{2i-1}=\lambda_{2i-1}+1$; if $\varepsilon(\lambda_{2i})=1$ and $ \lambda_{2i}>\lambda_{2i+1}$, then $\tilde{\lambda}_{2i}=\lambda_{2i}-1$; otherwise $\tilde{\lambda}_{i}={\lambda}_{i}$. Note that $\Psi_{SO(2n)}^2((\lambda,\varepsilon))=\tilde{\lambda}$ with $\tilde{\lambda}_i$ all even iff $\lambda=\tilde{\lambda}$ and $\varepsilon(\lambda_i)=0$ for all $\lambda_i$; for the two degenerate classes $\rc_1,\rc_2$ corresponding to such a $(\lambda,\varepsilon)$, we define $\Psi_G^2(\rc_i)=\tilde{\rc}_i$ by $\gamma_G^2(\rc_i)=\gamma_G^1(\tilde{\rc}_i)$, $i=1,2$.

Assume  that $G=SO(N)$ and $\rc=(\lambda,\chi)\in\Omega_{\Lg}^2$.  Let $k\geq 0$ be the unique integer such that $\lambda_{2k+2}=\lambda_{2k+1}-1$ (when $N$ is odd), and $k=\infty$ (when $N$ is even).
Define  $\Psi_\Lg^2(\rc)=\tilde{\lambda}$ as follows.
 If $\chi(\lambda_{1})>\frac{\lambda_{1}+1}{2}$, then $\tilde{\lambda}_1=2\chi(\lambda_1)-1$; if  $i\leq k$, $\chi(\lambda_{2i})>\frac{\lambda_{2i}+1}{2}$ and $\chi(\lambda_{2i})>\chi(\lambda_{2i+1})$, then $$\tilde{\lambda}_{2i}=\left\{\begin{array}{lll}\lambda_{2i}-\chi(\lambda_{2i})+\chi(\lambda_{2i+1})&\text{ if }\chi(\lambda_{2i})-\lambda_{2i}+\chi(\lambda_{2i+1})\geq2\\
2(\lambda_{2i}-\chi(\lambda_{2i}))+1&\text{ if }\chi(\lambda_{2i})-\lambda_{2i}+\chi(\lambda_{2i+1})\leq1\end{array}\right.;$$ if  $i\leq k$, $\chi(\lambda_{2i+1})>\frac{\lambda_{2i+1}+1}{2}$ and $\lambda_{2i+1}-\chi(\lambda_{2i+1})<\lambda_{2i}-\chi(\lambda_{2i})$, then $$\tilde{\lambda}_{2i+1}=\left\{\begin{array}{lll}\lambda_{2i}-\chi(\lambda_{2i})+\chi(\lambda_{2i+1})&\text{ if }\chi(\lambda_{2i})-\lambda_{2i}+\chi(\lambda_{2i+1})\geq2\\
2\chi(\lambda_{2i+1})-1&\text{ if }\chi(\lambda_{2i})-\lambda_{2i}+\chi(\lambda_{2i+1})\leq1\end{array}\right.;$$ if $i\geq k+1$, $\tilde{\lambda}_{2i}=\lambda_{2i+1}$;
 otherwise $\tilde{\lambda}_{i}={\lambda}_{i}$.
Note that $\Psi_{\mathfrak{so}(2n)}^2((\lambda,\chi))=\tilde{\lambda}$ with $\tilde{\lambda}_i$ all even iff $\lambda=\tilde{\lambda}$ and $\chi(\lambda_i)=\lambda_i/2$ for all $i$; for the two degenerate orbits $\rc_1,\rc_2$ corresponding to such a $(\lambda,\chi)$, we define $\Psi_\Lg^2(\rc_i)=\tilde{\rc}_i$ by $\gamma_\Lg^2(\rc_i)=\gamma_G^1(\tilde{\rc}_i)$, $i=1,2$.

We show  that

\qquad (a) {\em $\gamma_{G}^1\circ\Psi_G^p=\Phi_G\circ\gamma_{G}^p\ (\text{resp. }\gamma_{G}^1\circ\Psi_\Lg^2=\Phi_G\circ\gamma_{\Lg}^2).$}\\[10pt]
If $p\neq 2$, (a) is clear. Assume that $p=2$. We verify $(a)$ for $G=SO(2n+1)$. The other cases are entirely similar. Let $\rc\in\Omega_G^2$ (resp. $\Omega_\Lg^2$). Assume that $\gamma_G^2(\rc)=(\alpha,\beta)$ (resp. $\gamma_\Lg^2(\rc)=(\alpha,\beta)$), $\Phi_G((\alpha,\beta))=(\tilde{\alpha},\tilde{\beta})$ and $(\gamma_{G}^1)^{-1}((\tilde{\alpha},\tilde{\beta}))=\tilde{\lambda}$.  Using the definition of $\Phi_G$, one easily shows that $\tilde{\beta}_i=\tilde{\alpha}_i+2$ iff $\beta_i=\alpha_i+2$; $\tilde{\beta}_i=\tilde{\alpha}_{i+1}$ iff $\beta_i\leq\alpha_{i+1}$ and $\beta_i+\alpha_{i+1}$ is even. Using this and the description of the map $\gamma_{SO(2n+1)}^1$ in \ref{charnot2}, one easily verifies that we have
 $$\tilde{\lambda}_{2i-1}=\left\{\begin{array}{ll}\alpha_i+\beta_i&\text{ if }\alpha_i=\beta_i-2\\\alpha_i+\beta_{i-1}&\text{ if }\alpha_i\geq\beta_{i-1}\\2\alpha_i+1&\text{ if }\beta_i-1\leq\alpha_i\leq\beta_{i-1}-1\end{array}\right.,\ \tilde{\lambda}_{2i}=\left\{\begin{array}{ll}\alpha_i+\beta_i&\text{ if }\beta_i=\alpha_i+2\\\alpha_{i+1}+\beta_{i}&\text{ if }\beta_i\leq\alpha_{i+1}\\2\beta_i-1&\text{ if }\alpha_{i+1}+1\leq\beta_i\leq\alpha_{i}+1\end{array}\right..$$

Assume first that $\rc=(\lambda,\varepsilon)\in\Omega_G^2=\Omega_{Sp(2n)}^2$. By the description of the map $\gamma_{Sp(2n)}^2$ in \ref{sec-1}, we have $${\lambda}_{2i-1}=\left\{\begin{array}{ll}\alpha_i+\beta_i&\text{ if }\alpha_i\leq\beta_i-1\\\alpha_i+\beta_{i-1}&\text{ if }\alpha_i\geq\beta_{i-1}+1\\2\alpha_i&\text{ if }\beta_i\leq\alpha_i\leq\beta_{i-1}\end{array}\right.,\ {\lambda}_{2i}=\left\{\begin{array}{ll}\alpha_i+\beta_i&\text{ if }\beta_i\geq\alpha_i+1\\\alpha_{i+1}+\beta_{i}&\text{ if }\beta_i\leq\alpha_{i+1}-1\\2\beta_i&\text{ if }\alpha_{i+1}\leq\beta_i\leq\alpha_{i}\end{array}\right..$$ Thus $\tilde{\lambda}_{2i-1}=\lambda_{2i-1}+1$ iff $\beta_i\leq\alpha_i\leq\beta_{i-1}-1$ iff $\varepsilon(\lambda_{2i-1})=1$ and $\lambda_{2i-1}<\lambda_{2i-2}$; $\tilde{\lambda}_{2i}=\lambda_{2i}-1$ iff $\alpha_{i+1}+1\leq\beta_i\leq\alpha_i$ iff $\varepsilon(\lambda_{2i})=1$ and $\lambda_{2i}>\lambda_{2i+1}$; otherwise $\tilde{\lambda}_{i}=\lambda_{i}$. It follows from the definition of $\Psi_G^2$ that $\tilde{\lambda}=\Psi_G^2(\rc)$ and thus (a) follows.

Assume now that $(\lambda,\chi)\in\Omega_\Lg^2$. It is enough to show that $\tilde{\lambda}=\Psi_\Lg^2(\rc)$. By the description of the map $\gamma_{\mathfrak{so}(2n+1)}^2$ in \ref{nilpchar2}, we have  that  $\lambda_{2k+2}=\lambda_{2k+1}-1$ iff $\beta_k>0$ and $\beta_{k+1}=0$, thus $\tilde{\lambda}_{2i-1}=\alpha_i=\lambda_{2i-1}$ for all $i\geq k+2$, $\tilde{\lambda}_{2i}=\alpha_{i+1}=\lambda_{2i+1}$ for all $i\geq k+1$. It is easy to check that $\tilde{\lambda}_{2k+1}$ is as in the definition of $\Psi_\Lg^2(\rc)$. 
Assume now that  $i\leq k$. Then $\tilde{\lambda}_{2i-1}={\lambda}_{2i-1}$ iff $\alpha_i\leq\beta_i-1$ (namely $\chi(\lambda_{2i-1})\leq\frac{\lambda_{2i-1}+1}{2}$) or $\beta_{i-1}=\beta_i\leq\alpha_i$ (namely $\lambda_{2i-2}-\chi(\lambda_{2i-2})=\lambda_{2i-1}-\chi(\lambda_{2i-1})$ and $\chi(\lambda_{2i-1})>\frac{\lambda_{2i-1}+1}{2}$), $\tilde{\lambda}_{2i}={\lambda}_{2i}$ iff $\alpha_i\leq\beta_i-1$ (namely $\chi(\lambda_{2i})\leq\frac{\lambda_{2i}+1}{2}$) or $\alpha_{i+1}=\alpha_i\geq\beta_i$ (namely $\chi(\lambda_{2i})=\chi(\lambda_{2i+1})$ and $\chi(\lambda_{2i})>\frac{\lambda_{2i}+1}{2}$). Moreover, if $\lambda_{2i-2}-\chi(\lambda_{2i-2})>\lambda_{2i-1}-\chi(\lambda_{2i-1})$ and $\chi(\lambda_{2i-1})>\frac{\lambda_{2i-1}+1}{2}$, or $\chi(\lambda_{2i})>\chi(\lambda_{2i+1})$ and $\chi(\lambda_{2i})>\frac{\lambda_{2i}+1}{2}$, then it is easy to check that $\tilde{\lambda}_{2i-1}$ and $\tilde{\lambda}_{2i}$ are as in the definition of $\Psi_\Lg^2(\rc)$. This completes the verification of (a) for $SO(2n+1)$.

\begin{proposition}\label{prop-psi}
Two classes $\rc_1,\rc_2\in \Omega^p_{G}$ (resp. $\Omega^2_{\Lg}$) lie in the same unipotent (resp. nilpotent) piece as defined in \cite{Lu4} if and only if $\Psi_G^p(\rc_1)=\Psi_G^p(\rc_2)$ (resp. $\Psi_\Lg^2(\rc_1)=\Psi_\Lg^2(\rc_2)$).
\end{proposition}
The proposition is clear when $p\neq 2$. If $G=Sp(2n)$, the proposition follows from \cite{Lu4} (see Lemma \ref{lem-sp1}) and the definition of $\Psi_G^p$ (resp. $\Psi_\Lg^2$). The proof of the proposition in the case where $p=2$ and $G=SO(N)$ is given in subsections \ref{pf-char2}-\ref{pf-char2-d}. Note that the proposition computes the pieces in classical groups explicitly. Another computation of the unipotent pieces is given in \cite{Lu8}. Now in view of  (a) and \ref{ssec-fiber1} (a), \ref{ssec-comb1} (b) follows from Proposition \ref{prop-psi}.

\subsection{}\label{pf-char2}
Assume that $p=2$ and $G=SO_Q(V)$ in the remainder of this section.

Let $\rc=(\lambda,\varepsilon)\in\Omega_G^2$ (resp. $\rc=(\lambda,\chi)\in\Omega_\Lg^2$).  Assume
$$\Psi_G^2(\rc)=\tilde{\rc}=\tilde{\lambda}\text{ (resp. }\Psi_\Lg^2(\rc)=\tilde{\rc}=\tilde{\lambda})$$where $\tilde{\lambda}=(\tilde{\lambda}_1\geq\tilde{\lambda}_2\geq\cdots)$. Suppose that $\Upsilon_{\rc}=(f_a)_{a\in\mathbb{N}}$ and $\Upsilon_{\tilde{\rc}}=(\tilde{f}_a)_{a\in\mathbb{N}}$ (see \ref{lem-o1}).
We show that

\qquad (a) {\em $f_a=\tilde{f}_a$ for all ${a\in\mathbb{N}}$}.\\[10pt]
Then Proposition \ref{prop-psi} follows from (a) and Lemma
\ref{lem-o1} (note that for $\tilde{\rc}_i=\tilde{\lambda}^i\in\Omega_G^1$ with $\Upsilon_{\tilde{\rc}_1}=\Upsilon_{\tilde{\rc}_2}$, by \ref{lem-o1} (a$'$)  we have $\tilde{\lambda}^1=\tilde{\lambda}^2$; note also that  $\rc$ is a degenerate class (if $G=SO(2n)$)   iff $\tilde{\rc}$ is a degenerate class iff $\tilde{f}_0=0$ and then by \cite{Lu4}  each degenerate class itself forms one piece).

We prove (a) by induction
on $\dim V$. Let $u\in\rc$ (resp. $x\in\rc$) and let $T=u-1$ (resp. $x$).
If $T=0$, (a) is obvious. Assume  from now on that $T\neq0$.
Let $V_*=(V_{\geq a})$,  $V'$, $u'$ (resp. $x'$) be associated to $u$ (resp. $x$)  and $m$, $e$, $f$ defined for $T$ as in
\ref{lem-o1}.   Let $\rc'$ be the class of $u'$ (resp. $x'$) in $G'=SO_{Q'}(V')$ (resp. $\Lg'=\mathfrak{so}_{Q'}(V')$) and let
$$\tilde{\rc}'=\Psi_{G'}^2(\rc')=\tilde{\lambda}'\text{ (resp. }\tilde{\rc}'=\Psi_{\Lg'}^2(\rc')=\tilde{\lambda}').$$
Suppose that $\Upsilon_{\rc'}=(f_a')$ and $\Upsilon_{\tilde{\rc}'}=(\tilde{f}_a').$  Since $\dim V'<\dim V$, by induction hypothesis
$f_a'=\tilde{f}_a'$ for all $a\in\mathbb{N}$. By the definition of $V_*$ we have that for all $a\in[0,m-1]$, $f_a=f_a'$  and thus $f_a=\tilde{f}'_a$.
We show  that

\qquad (b) {\em $\begin{array}{l}
\tilde{\lambda}_1=m+1, m_{\tilde{\lambda}}(\tilde{\lambda}_1)=f_m,  m_{\tilde{\lambda}'}(\tilde{\lambda}_1)=0,\
m_{\tilde{\lambda}'}(\tilde{\lambda}_1-2)=m_{\tilde{\lambda}}(\tilde{\lambda}_1-2)+m_{\tilde{\lambda}}(\tilde{\lambda}_1),
\\[5pt]m_{\tilde{\lambda}'}(i)=m_{\tilde{\lambda}}(i) \text{ for all }i\neq \tilde{\lambda}_1,\tilde{\lambda}_1-2.\end{array}$}\\[10pt]
It then follows from (b) and \ref{lem-o1} (a$'$) that $\tilde{f}_a=0$ for all $a\geq m+1$, $\tilde{f}_m=f_m$, and  that $\tilde{f}_a=\tilde{f}_a'$ for all $a\in[0,m-1]$.  Hence (a) follows (note that $f_a$=0 for all $a\geq m+1$).

The proof of (b) for $\rc\in\Omega_G^2$ is given in subsections \ref{pf-char2-3}-\ref{pf-char2-4} and that for $\rc\in\Omega_\Lg^2$ is given in subsections \ref{pf-char2-2}-\ref{pf-char2-d}. For a subspace $W\subset V$ and a subset $E\subset W$, we denote $$W^\p=\{v\in V|\la v,W\ra=0\},\ E^{\p_W}=\{v\in W|\la v,E\ra=0\}.$$

\subsection{} \label{pf-char2-3}
Assume that $\rc=(\lambda,\varepsilon)\in\Omega_G^2$ in this subsection.  We have $e=\lambda_1$; $e=2f-2$ if $\varepsilon(e)=1$, $e=2f-1$ if $\varepsilon(e)=\omega$, and $e=2f$ if $\varepsilon(e)=0$. Note that $e\geq 2$ since $T\neq 0$. We keep the notations in \ref{pf-char2}. We can compute $\rc'=(\lambda',\varepsilon')$ in various cases as follows. 

\noindent(i) $e=2f-2$. We have\\[5pt]
\indent{\em if $m_\lambda(e)$ is even, then $m_{\lambda'}(e)=m_\lambda(e)-2,\ m_{\lambda'}(e-1)=m_{\lambda}(e-1)+2,\  m_{\lambda'}(i)=m_{\lambda}(i) \text{ for }i\notin\{ e,e-1\},\
\varepsilon'(e)\leq0,\ \varepsilon'(\lambda_i)=\varepsilon(\lambda_i)\text{ for }\lambda_i\neq e$};\\[5pt]
\indent{\em if $m_\lambda(e)$  is odd, then $m_{\lambda'}(e)=m_\lambda(e)-1,\ m_{\lambda'}(e-2)=m_{\lambda}(e-2)+1\ (\text{if }e>2),\  m_{\lambda'}(i)=m_{\lambda}(i)$ for $i\notin\{ e,e-2\},\ 
\varepsilon'(e)\leq 0,\ \varepsilon'(e-2)=1\ (\text{if }e>2),\ \varepsilon'(\lambda_i)=\varepsilon(\lambda_i)\text{ for }\lambda_i\notin\{ e,e-2\}.$}

\noindent(ii) $e=2f-1$ and $\varepsilon(e-1)=1$. We have\\[5pt]
 \indent{\em if $m_\lambda(e-1)$ is even, then $\ m_{\lambda'}(e)=0,\ m_{\lambda'}(e-1)=m_{\lambda}(e-1)-2,\ m_{\lambda'}(e-2)=m_{\lambda}(e-2)+m_{\lambda}(e)+2\ (\text{if }e>2), \ m_{\lambda'}(i)=m_{\lambda}(i) \text{ for }i\notin\{ e,e-1,e-2\},\
\varepsilon'(e-1)\leq 0,\ \varepsilon'(\lambda_i)=\varepsilon(\lambda_i)\text{ for }\lambda_i\neq e-1$;}\\[5pt]
\indent{\em if $m_\lambda(e-1)$  is odd, then $m_{\lambda'}(e)=0,\ m_{\lambda'}(e-1)=m_{\lambda}(e-1)-1,\ m_{\lambda'}(e-2)=m_{\lambda}(e-2)+m_{\lambda}(e)\ (\text{if }e>2),\ m_{\lambda'}(e-3)=m_{\lambda}(e-3)+1\ (\text{if }e>3),\ m_{\lambda'}(i)=m_{\lambda}(i) \text{ for }i\notin\{ e,e-1,e-2,e-3\},\
\varepsilon'(e-1)\leq 0,\ \varepsilon'(e-3)=1\ (\text{if }e>3),\ \varepsilon'(\lambda_i)=\varepsilon(\lambda_i)\text{ for }\lambda_i\notin\{e-1,e-3\}.$}

\noindent(iii)  $e=2f-1$ and $\varepsilon(e-1)\leq0$, or $e=2f$. We have\\[5pt]
\indent {\em $m_{\lambda'}(e)=0,\ m_{\lambda'}(e-2)=m_{\lambda}(e-2)+m_{\lambda}(e)\ (\text{if }e>2),\ m_{\lambda'}(i)=m_{\lambda}(i) \text{ for }i\notin\{ e,e-2\}; \  \varepsilon'(e-2)=1\text{ iff } \varepsilon(e-2)=1,\
\varepsilon'(\lambda_i)=\varepsilon(\lambda_i)\text{ for }\lambda_i\neq e-2.$}

\noindent In fact  the result in case (i) follows from \cite[2.4 (ii) (iii)]{Lu3} and that in case (iii) follows from \cite[2.4 (i)]{Lu3} (in these cases our $V'$ is the same as $V'$ in \cite{Lu3}).
Assume now that we are in case (ii).  Let $T'=u'-1$ be the map on $V'$ induced by $T$. We have a decomposition of $V$ into mutually orthogonal $T$-stable subspaces $V=W_0\oplus W_1$ such that
$$\lambda_{T_0}=e^{m_\lambda(e)},\ T^{e-1}W_1=0,\ \varepsilon_{T_1}(e-1)=1,$$
 where $T_i=T|_{W_i}=(\lambda_{T_i},\varepsilon_{T_i})$, $i=0,1$. Note that $Q(T^{f-1}w)=0$ for all $w\in\ker T^{e-1}\cap{W_0}$ and $Q(T^{e-1}W_0)=0$ (as $e-1\geq f$). Let $K_{W_1}=\{v\in W_1|Q(T^{f-1}v)=0\}$ and $L_{W_1}=K_{W_1}^{\p_{W_1}}\cap Q^{-1}(0)$. Then $$ V_{\geq -m+1}=(\ker T^{e-1}\cap{W_0})\oplus K_{W_1},\ V_{\geq m}=T^{e-1}W_0\oplus L_{W_1}.$$  Thus we have a natural decomposition of $V'$ into mutually orthogonal $T'$-stable subspaces $$V'=W_0'\oplus W_1', \ W_0'=(\ker T^{e-1}\cap{W_0})/T^{e-1}W_0,\ W_1'=K_{W_1}/L_{W_1}.$$   It is easy to see that $\lambda_{T'|_{W_0'}}=(e-2)^{m_\lambda(e)}.$  Let $e_{T_1}$, $f_{T_1}$ be defined for $T_1$ as $e,f$ for $T$. We have $e_{T_1}=e-1$ and $f_{T_1}=f$. Thus we can apply the result in case (i) to $T_1$ on $W_1$ and then the result in case (ii) follows.

\subsection{}\label{pf-char2-4}
We keep the notations in \ref{pf-char2-3}.
Using the definition of $\Psi_G^2$ and the description of $\rc'$ in \ref{pf-char2-3}, we can compute $\tilde{\lambda}$ and $\tilde{\lambda}'$ in each case (i)-(iii) as follows. Let $d_i=m_\lambda(e-i)$.

Assume first that we are in case (i) with $d_0$ odd (resp. in case (ii) with $d_1$ odd). Let $d=d_0+d_1$ (resp. $d=d_0+d_1+d_2$). Then $d$ is odd. Note that $\lambda_i'=\lambda_i$ for all $i\geq d+1$. For $2i-1\geq d+2$, $\lambda_{2i-1}'<\lambda_{2i-2}'(\leq \lambda_d')$ and $\varepsilon'(\lambda_{2i-1}')=1$ iff $\lambda_{2i-1}<\lambda_{2i-2}$ and $\varepsilon(\lambda_{2i-1})=1$, thus $\tilde{\lambda}_{2i-1}=\tilde{\lambda}_{2i-1}'$; for $2i\geq d+1$, $\lambda_{2i}'>\lambda_{2i+1}'$ and $\varepsilon'(\lambda_{2i}')=1$ iff $\lambda_{2i}>\lambda_{2i+1}$ and $\varepsilon(\lambda_{2i})=1$ (note that if $\lambda_{2i}'=\lambda_d'$, then $m_\lambda(\lambda_d')$ is odd and thus $\varepsilon(\lambda_d')=1$),  thus $\tilde{\lambda}_{2i}=\tilde{\lambda}_{2i}'$. We have shown that $\tilde{\lambda}_i=\tilde{\lambda}_i'\text{ for all }i\geq d+1.$ Let $\tilde{\lambda}^1=(\tilde{\lambda}_1,\ldots,\tilde{\lambda}_{d})$ and $\tilde{\lambda}^{1'}=(\tilde{\lambda}'_1,\ldots,\tilde{\lambda}'_{d}).$ We have

\quad in case (i) (with $d_0$ odd) {\em $\tilde{\lambda}^1=(e+1)e^{d_0-1}(e-1)^{d_1}$, $\tilde{\lambda}^{1'}=e^{d_0-1}(e-1)^{d_1+1}$;}\\[5pt]
\indent \quad in case (ii) (with $d_1$ odd) {\em $\tilde{\lambda}^1=e^{d_0+1}(e-1)^{d_1-1}(e-2)^{d_2}$, $\tilde{\lambda}^{1'}=(e-1)^{d_1-1}(e-2)^{d_0+d_2+1}$.}

 Assume now that we are in the remaining cases.  Let $d=d_0+d_1$. Then $d$ is even. Note that for all $i\geq d+1$, $\lambda_i'=\lambda_i$, $\varepsilon'(\lambda_i')=1$ iff $\varepsilon(\lambda_i)=1$. Hence $\tilde{\lambda}_i=\tilde{\lambda}_i'\text{ for all }i\geq d+2.$ Let $\tilde{\lambda}^1=(\tilde{\lambda}_1,\ldots,\tilde{\lambda}_{d+1}),\  \tilde{\lambda}^{1'}=(\tilde{\lambda}'_1,\ldots,\tilde{\lambda}'_{d+1}).$ If we are not in case (iii) with $\varepsilon(e-2)=1$, then $\tilde{\lambda}_{d+1}=\tilde{\lambda}'_{d+1}$ since $\varepsilon'({\lambda}'_{d+1})=1$ implies that ${\lambda}_{d+1}'<\lambda_d'$ and thus $\varepsilon({\lambda}_{d+1})=1$ and $\lambda_{d+1}<\lambda_d$. We have

\ in case (i) (with $d_0$ even) {\em $
\tilde{\lambda}^1=(e+1)e^{d_0-2}(e-1)^{d_1+1}\tilde{\lambda}_{d+1}, \  \tilde{\lambda}^{1'}=e^{d_0-2}(e-1)^{d_1+2}\tilde{\lambda}_{d+1}$;}

\ in case (ii) (with $d_1$ even) {\em $
\tilde{\lambda}^1=e^{d_0+1}(e-1)^{d_1-2}(e-2)\tilde{\lambda}_{d+1}\ \tilde{\lambda}^{1'}=(e-1)^{d_1-2}(e-2)^{d_0+2}\tilde{\lambda}_{d+1}
$;}

\ in case (iii) {\em $\left\{\begin{array}{ll} 
\tilde{\lambda}^1=e^{d_0}(e-1)^{d_1+1},\ 
\tilde{\lambda}^{1'}=(e-1)^{d_1+1}(e-2)^{d_0}&\text{if }\varepsilon(e-2)=1\\
\tilde{\lambda}^1=e^{d_0}(e-1)^{d_1}\tilde{\lambda}_{d+1},\ 
\tilde{\lambda}^{1'}=(e-1)^{d_1}(e-2)^{d_0}\tilde{\lambda}_{d+1}&\text{if }\varepsilon(e-2)\leq 0.\end{array}\right.$}

We have $m=e$ in case (i) and $m=e-1$ in case (ii) (iii). Now \ref{pf-char2} (b) for $\rc\in\Omega_G^2$ follows from the above description of $\tilde{\lambda},\tilde{\lambda}'$ and that (see \cite[2.6 (i), (iii),  (vi)]{Lu3} )

\qquad {\em $
f_m=\left\{\begin{array}{ll}1&\text{ if }e=2f-2,\\
m_{{\lambda}}(e)+1&\text{if }e=2f-1\text{ and }\varepsilon(e-1)=1,\\m_{{\lambda}}(e)&\text{if }e=2f-1\text{ and }\varepsilon(e-1)\leq 0,\text{ or }e=2f\\\end{array}\right.$.}

\subsection{} \label{pf-char2-2} Assume $\rc=(\lambda,\chi)\in\Omega_\Lg^2$ in the remainder of this section. We keep the notations in \ref{pf-char2}. Let $\chi_T:\mathbb{N}\to\mathbb{N}$ be the function defined for $T$ as follows $$ \chi_T(a)=\min\{s\in\mathbb{N}|T^{a}v=0\Rightarrow Q(T^{s}v)=0,v\in V\}.$$ Then  $\chi(\lambda_i)=\chi_T(\lambda_i)$; we write $\chi_T=\chi$. We have

\qquad (a) $e=\lambda_1,\ f=\chi(\lambda_1).$\\[10pt]
 We show in this subsection that

\qquad (b) {\em $
f_m=\left\{\begin{array}{ll}1&\text{ if }e<2f-1,\\
m_{{\lambda}}(e)+1&\text{if }e=2f-1\text{ and }\chi(e-1)=f,\\
m_{{\lambda}}(e)&\text{if }e=2f-1\text{ and }\chi(e-1)=f-1,\text{ or } e=2f.\end{array}\right.$}\\[10pt]
Recall that $f_m=\dim V_{\geq m} $ and $V_{\geq m}= V_{\geq -m+1}^\p\cap Q^{-1}(0)$ (see \ref{lem-o1}). Note that $e\geq 2$ since $T\neq 0$.

 Suppose  that $e<2f-1$. Consider the map $\rho:V\to\tk$, $v\mapsto \sqrt{Q(T^{f-1}v)}$, where $\sqrt{\ }$ is a chosen square root on $\tk$. It is easy to show that $\rho$ is linear. Thus if $R=0$, then $V_{\geq m}=(\ker\rho)^\p$ is a line; if $R\neq 0$, $(\ker\rho)^\p$ is a two dimensional subspace of $V$ containing $R$ and $(\ker\rho)^\p=V_{\geq m}\oplus R$. In each case, we have $\dim V_{\geq m}=1$.
 
 Suppose that $e=2f-1 \text{ and }\chi(e-1)=f$. Let $E$ be a complement to $\ker T^{e-1}$ in $V$ and $W=E+TE+\cdots+T^{e-1}E$. Similar  argument as in \cite[1.9]{Lu3} shows that $\la,\ra|_W$ is non-degenerate and we have a decomposition $V=W\oplus Y$, where $Y=W^\p$ is $T$-stable and $T^{e-1}Y=0$. It is then easy to show that $V_{\geq -m+1}=(\ker T^{e-1}\cap W)\oplus\{v\in Y|Q(T^{f-1}v)=0\}$ and thus $V_{\geq m}=T^{e-1}V\oplus L_Y$, where $L_Y\subset Y$ is a line (we apply the discussion in the first case for $Y$). Hence $f_m=m_\lambda(e)+1$.

 Suppose that $e=2f-1 \text{ and }\chi(e-1)=f-1$, or $e=2f$. Then $V_{\geq -m+1}=\ker T^{e-1}$. Note that $T^{e-1}V\cap R=0$ since $Q|_{T^{e-1}V}=0$ (as $e-1\geq f$). Thus $V_{\geq -m+1}^\p=T^{e-1}V\oplus R$ and $V_{\geq m}=T^{e-1}V$. It follows that $f_m=m_\lambda(e)$.

\subsection{}\label{pf-char2-5}We keep the notations in \ref{pf-char2-2}.  Let  $j\geq 0$ be the unique integer such that  $$\chi(e-j)=f\text{ and }\chi(e-j-1)<f.$$We describe $\rc'=(\lambda',\chi')$  in various cases as follows.

\noindent(i) $e<2f-1$. We have\\[5pt]
 \indent {\em $\chi'(e-k)=f-1\text{ for }k\in[0,j],\ \chi'(\lambda_i)=\chi(\lambda_i)\text{ for }\lambda_i\leq e-j-1$;}\\[5pt]
 \indent{\em if $m_\lambda(e-j)$ is even, then $
m_{\lambda'}(e-j)=m_\lambda(e-j)-2,\ m_{\lambda'}(e-j-1)=m_{\lambda}(e-j-1)+2$ (if $e>j+1$), $ m_{\lambda'}(i)=m_{\lambda}(i)\text{ for }i\notin\{ e-j,e-j-1\},$}\\[5pt]
\indent{\em 
if $m_\lambda(e-j)$ is odd, then $m_{\lambda'}(e-j)=m_\lambda(e-j)-1,\ m_{\lambda'}(e-j-2)=m_{\lambda}(e-j-2)+1$ (if $e>j+2$), $m_{\lambda'}(i)=m_{\lambda}(i) \text{ for }i\notin\{ e-j,e-j-2\}.$}

\noindent(ii) $e=2f-1$ and $\chi(e-1)=f$. We have \\[5pt]
\indent{\em $
\chi'(e-k)=f-1\text{ for }k\in[0,j],\ \chi'(\lambda_i)=\chi(\lambda_i)\text{ for }\lambda_i\leq e-j-1$;}\\[5pt]
\indent{\em if $m_\lambda(e-j)$  is even, then $m_{\lambda'}(e)=0,\ m_{\lambda'}(e-2)=m_{\lambda}(e-2)+m_{\lambda}(e)+2\delta_{j,1}-2\delta_{j,2}$ (if $e>2$), $ m_{\lambda'}(e-j)=m_{\lambda}(e-j)-2+\delta_{j,2}m_\lambda(e),\ m_{\lambda'}(e-j-1)=m_{\lambda}(e-j-1)+2+\delta_{j,1}m_\lambda(e)$ (if $e>j+1$), $ m_{\lambda'}(i)=m_{\lambda}(i)\text{ for }i\notin\{ e,e-2,e-j,e-j-1\},$}\\[5pt]
\indent{\em 
if $m_\lambda(e-j)$  is odd, then $m_{\lambda'}(e)=0,\ m_{\lambda'}(e-2)=m_{\lambda}(e-2)+m_{\lambda}(e)-\delta_{j,2}$ (if $e>2$), $  m_{\lambda'}(e-j)=m_{\lambda}(e-j)-1+m_{\lambda}(e)\delta_{j,2},\ m_{\lambda'}(e-j-2)=m_{\lambda}(e-j-2)+1 \ (\text{if }e>j+2),\ m_{\lambda'}(i)=m_{\lambda}(i) \text{ for }i\notin\{ e,e-2,e-j,e-j-2\}$.}

\noindent (iii) $e=2f-1$ and $\chi(e-1)=f-1$, or $e=2f$. We have\\[5pt]
\indent{\em $
\chi'(\lambda_i)=\chi(\lambda_i)\text{ for }\lambda_i\leq e-1;$}\\[5pt]\indent {\em $m_{\lambda'}(e)=0,\ m_{\lambda'}(e-2)=m_{\lambda}(e-2)+m_{\lambda}(e)\ (\text{if }e>2),\ m_{\lambda'}(i)=m_{\lambda}(i) \text{ for }i\notin\{ e,e-2\}.$}

We explain the computation of $\rc'$ in more detail in the remainder of this subsection. 
Recall from \cite[3.7]{Hes} that we have a decomposition of $V$ into mutually orthogonal $T$-stable subspaces $V=W(1)\oplus W(2)\oplus\cdots\oplus W(r)$ such that $m_\lambda(i)=\sum _{a\in[1,r]}m_{\lambda^a}(i)$, $\chi(i)=\max_a\chi_a(i)$, where $T|_{W(a)}=(\lambda^a,\chi_a)$. Moreover, $T|_{W(a)}=W_{\chi(\lambda_i)}(\lambda_i)$ (for some $i$) for $a\in[1,r-1]$,  $T|_{W(r)}=W_{\chi(\lambda_i)}(\lambda_i)$ (for some $i$) if $R=0$, and $T|_{W(r)}\cong D(\lambda_{2k+1})$  if $R\neq0$,
where

 (a$_1$) $T|_W=W_{l}(s)$ means that  there exist $v_0,w_0\in W$ such that $W=\text{span}\{T^iv_0,T^{i}w_0,i\in[0,s-1]\}, \ \la T^iv_0,w_0\ra=\delta_{i,s-1},\  Q(T^iv_0)=\delta_{i,l-1},\ Q(T^iw_0)=0$; we have $\chi_{T|_W}(i)=\max(0,\min(i-s+l,l))$;

 (a$_2$) $T|_W=D(s)$ means that there exist $v_0,w_0\in W$ such that $W=\text{span}\{T^iv_0, T^iw_0,i\in[0,s-2],T^{s-1}v_0\},\ \la T^iv_0,w_0\ra=\delta_{i,s-2},\  Q(T^iv_0)=\delta_{i,s-1},\ Q(T^iw_0)=0$; we have $\chi_{T|_W}(i)=\min(i,s)$.
 
 The following facts will be used  in the  computation of $\rc'$.

(b$_1$) Let $W$ be a $T$-stable subspace of $V$ such that $T|_W=W_f(e-j)$ with  $f>\frac{e-j+1}{2}$ (resp.  $T|_W=D(e-j)$ with $f=e-j$). Let  $K_W=\{v\in W|Q(T^{f-1}v)=0\}$, $L_W=K_W^{\p_W}\cap Q^{-1}(0)$, $W'=K_W/L_W$ and let $T_1'$ be the map on $W'$ induced by $T$. Using the basis for $W$ chosen as in (a$_1$) (resp. (a$_2$)), one can easily check that  $$T_1'|_{W'}=W_{f-1}(e-j-1)\ (\text{resp. }T_1'|_{W'}=D(e-j-1)).$$

(b$_2$) Let $W$ be a $T$-stable subspace of $V$ such that $T|_W=W_f(e)^{a}$ (an orthogonal decomposition into $a$ copies of $W_f(e)$), where $f\leq\frac{e+1}{2}$. Let $W'=(\ker T^{e-1}\cap W)/T^{e-1}W$ and $T_1'$ be the map on $W'$ induced by $T$. Using the basis for $W$ chosen as in (a$_1$), one can easily check that  $$T_1'|_{W'}=W_{f-1}(e-2)^{a}.$$
\indent(c) If $W_1$ and $W_2$ are two $T$-stable subspaces of $V$ and $T|_{W_i}=(\lambda^i,\chi^i)$ with $\lambda^1=\lambda^2$ and $\chi^1<\chi^2\leq \chi$, then $V=W_1\oplus W_1^\p\cong W_2\oplus W_1^\p$ (see \cite[Lemma 3.6]{Hes}).

Now we are ready to compute $\rc'$ in various cases. Let $d_i=m_\lambda(e-i)$. We have $d_j>0$ and $d_i, i\in[0,j-1]$ is even (since $\chi(e-i)=f\leq e-j<e-i$). Note that if $d_j$ is odd, then $\dim V$ is odd and $e-j=f$.

Assume first that we are in case (i) and that $d_j$ is even (resp. odd).  We have a decomposition of $V$ into mutually orthogonal $T$-stable subspaces $V=W\oplus Y$  such that (we use (c))

 \qquad $T|_W=W_{f}(e-j)\ (\text{resp. } T|_W=D(e-j))$  and $Q(T^{f-1}Y)=0$.\\[10pt]
   We have $V_{\geq -m+1}=K_W\oplus Y$ and $ V_{\geq m}=L_W,$ where $K_W=\{v\in W|Q(T^{f-1}v)=0\}$ and $L_W=K_W^{\p_W}\cap Q^{-1}(0)$. Hence we have a natural decomposition of $V'$ into mutually orthogonal $T'$-stable subspaces $V'=W'\oplus Y$, where $W'=K_W/L_W$. Moreover (see (b$_1$))

 \qquad   $T'|_{W'}=W_{f-1}(e-j-1)$ (resp. $D(e-j-1)$), $T'|_Y=T|_Y$. \\[10pt]
We have  $\chi'(i)=\max(\chi_{{T'}|_{W'}}(i),\chi_{T|_Y}(i))$ and $\chi(i)=\max(\chi_{T|_W}(i),\chi_{T|_Y}(i))$. If $i\leq e-j-1$, then  $\chi_{{T'}|_{W'}}(i)=\chi_{T|_W}(i)$ and thus $\chi'(i)=\chi(i)$. Now for $k\in[0,j]$, $\chi_{T'|_{W'}}(e-k)=f-1$ and thus $\chi'(e-k)=f-1$   (note that $\chi_{T|_Y}(e-k)\leq f-1$).

Assume that we are in case (ii).  Then $j\geq 1$. Assume that $d_j$ is even (resp. odd).  We have a decomposition of $V$ into mutually orthogonal $T$-stable subspaces $V=W_0\oplus W_1\oplus Y$ such that (we use (c))

$T|_{W_0}=W_{f}(e)^{\frac{d_0}{2}} $, $T_{W_1}=W_f(e-j)$ (resp. $D(e-j)$),  and $Q(T^{f-1}Y)=0$.\\[10pt]
We have $V_{\geq -m+1}=(\ker T^{e-1}\cap{W_0})\oplus K_{W_1}\oplus Y$ and $V_{\geq m}=T^{e-1}W_0\oplus L_{W_1},$ where $K_{W_1}=\{v\in W_1|Q(T^{f-1}v)=0\}$ and $L_{W_1}=K_{W_1}^{\p_{W_1}}\cap Q^{-1}(0)$. Hence we have a natural decomposition of $V'$ into mutually orthogonal $T'$-stable subspaces $V'=W_0'\oplus W_1'\oplus Y$, where $W_0'=(\ker T^{e-1}\cap{W_0})/T^{e-1}W_0,\ W_1'=K_{W_1}/L_{W_1}$, and (see (b$_1$) and (b$_2$))

\qquad$T'|_{W_0'}=W_{f-1}(e-2)^{\frac{d_0}{2}} ,\ T'|_{W_1'}=W_{f-1}(e-j-1)\text{ (resp. }D(e-j-1)),\ T'|_Y=T|_Y$.\\[10pt]
We have $\chi'(i)=\max(\chi_{T'|_{W_0'}}(i),\chi_{T'|_{W_1'}}(i),\chi_{T|_Y}(i))$ and $\chi(i)=\max(\chi_{T|_{W_0}}(i),\chi_{T|_{W_1}}(i),\chi_{T|_Y}(i))$. Note that $\chi_{T|_Y}(i)\leq f-1$. For $i\leq e-j-1$, $\chi_{T|_{W_0}}(i)\leq\chi_{T'|_{W_0'}}(i)\leq\chi_{T'|_{W_1'}}(i)=\chi_{T|_{W_1}}(i)(=\max(0,i-e+f+j))$ (since $j\geq 1$) and thus $\chi'(i)=\chi(i)$; for $e-j\leq i\leq e-2$, $\chi'(i)=\max(i-e+f+1,f-1)=f-1$; $\chi'(e-1)=f-1$.

Assume now that we are in case (iii).  We have a decomposition of $V$ into mutually orthogonal $T$-stable subspaces $V=W\oplus Y$  such that

\qquad $T|_W=W_f(e)^{\frac{d_0}{2}} $ and $T^{e-1}Y=0$.\\[10pt] We have $V_{\geq -m+1}=(\ker T^{e-1}\cap W)\oplus Y$ and $V_{\geq m}=T^{e-1}W$. Hence we have a natural decomposition of $V'$ into mutually orthogonal $T'$-stable subspaces $V'=W'\oplus Y$, where $W'=(\ker T^{e-1}\cap W)/T^{e-1}W$, and  (see (b$_2$))

\qquad $T'|_{W'}=W_{f-1}(e-2)^{\frac{d_0}{2}} ,\ T'|_Y=T|_Y$. \\[10pt]
We have  $\chi'(i)=\max(\chi_{{T'}|_{W'}}(i),\chi_{T|_Y}(i))$ and $\chi(i)=\max(\chi_{T|_W}(i),\chi_{T|_Y}(i))$. For $0<\lambda_i\leq e-2$ (if $e=2f-1$), or $0<\lambda_i\leq e-1$ (if $e=2f$), we have $\lambda_i-\chi(\lambda_i)\leq\frac{\lambda_i}{2}<e-f$ and thus $\chi_{T|_W}(\lambda_i)=\max(0,\lambda_i-e+f)<\chi(\lambda_i)$, which implies that $\chi_{T|_Y}(\lambda_i)=\chi(\lambda_i)$ and thus $\chi'(\lambda_i)=\max(\max(\lambda_i-e+f+1,0),\chi_{T|_Y}(\lambda_i))=\chi(\lambda_i)$.  Now if $\lambda_i=e-1$ and $e=2f-1$, then $\chi'(\lambda_i)=\max(f-1,\chi_{T|_Y}(\lambda_i))=\chi(\lambda_i)$ (we have $\chi_{T|_Y}(\lambda_i)\leq f-1$ and $\chi(e-1)=f-1$).

\subsection{}\label{pf-char2-d}We keep the notations in \ref{pf-char2-5}. Using the definition of $\Psi_\Lg^2$ and the  description of $\rc'$ in \ref{pf-char2-5}, we compute $\tilde{\lambda}$ and $\tilde{\lambda}'$ in each case (i)-(iii) as follows.

Assume first that we are in case (i) or (ii) with $d_j$ even. Then $d_a$ is even for $a\in[0,j]$. Let $d=\sum_{a\in[0,j]}d_a$. We have $\tilde{\lambda}_i=\tilde{\lambda}_i'$ for all $i\geq d$, since $\lambda_i=\lambda_i'\leq e-j-1$ and $\chi(\lambda_i)=\chi(\lambda_i')$ for all $i\geq d+1$ and $\lambda_d-\chi(\lambda_d)=\lambda_d'-\chi'(\lambda_d')$. Let $\tilde{\lambda}^1=(\tilde{\lambda}_1,\ldots,\tilde{\lambda}_{d-1})$ and $\tilde{\lambda}^{1'}=(\tilde{\lambda}'_1,\ldots,\tilde{\lambda}'_{d-1})$. We have

\quad in case (i) ($d_j$ even) $\left\{\begin{array}{l}
\tilde{\lambda}^1=(2f-1)e^{d_0}(e-1)^{d_1}\cdots(e-j+1)^{d_{j-1}}(e-j)^{d_{j}-2}\\
\tilde{\lambda}^{1'}=(2f-3)e^{d_0}(e-1)^{d_1}\cdots(e-j+1)^{d_{j-1}}(e-j)^{d_{j}-2}
\end{array}\right.;$

\quad in case (ii) ($d_j$ even) $\left\{\begin{array}{ll}
\tilde{\lambda}^1=e^{d_0+1}(e-1)^{d_1}\cdots(e-j+1)^{d_{j-1}}(e-j)^{d_j-2}\\\tilde{\lambda}^{1'}=(e-2)^{d_0+1}(e-1)^{d_1}\cdots(e-j+1)^{d_{j-1}}(e-j)^{d_j-2}.\end{array}\right.$

Assume that we are in case (i) or (ii) with $d_j$ odd. Then $d_{j+1}$ is odd. We have $2k+1=\sum_{a\in[0,j]}d_a$. Let $k'$ be the unique integer such that $\lambda_{2k'+2}'=\lambda_{2k'+1}'+1$. We have $2k'+1=\sum_{a\in[0,j+1]}d_a-1$ and ${\lambda}_i={\lambda}_i'$ for all $i\geq 2k'+2$ and thus
$\tilde{\lambda}_i=\tilde{\lambda}_i'$ for all $i\geq 2k'+2$. Let $\tilde{\lambda}^1=(\tilde{\lambda}_1,\ldots,\tilde{\lambda}_{2k'+1})$ and $\tilde{\lambda}^{1'}=(\tilde{\lambda}'_1,\ldots,\tilde{\lambda}'_{2k'+1})$. We have 

\quad in case (i) ($d_j$ odd) $\left\{
\begin{array}{l}
\tilde{\lambda}^1=(2f-1)e^{d_0}(e-1)^{d_1}\cdots(e-j+1)^{d_{j-1}}(e-j)^{d_j-1}(e-j-1)^{d_{j+1}-1},\\
\tilde{\lambda}^{1'}=(2f-3)e^{d_0}(e-1)^{d_1}\cdots(e-j+1)^{d_{j-1}}(e-j)^{d_j-1}(e-j-1)^{d_{j+1}-1}\end{array}\right.$;

\quad in case (ii) ($d_j$ odd) $\left\{\begin{array}{l}
\tilde{\lambda}^1=e^{d_0+1}(e-1)^{d_1}\cdots(e-j+1)^{d_{j-1}}(e-j)^{d_j-1}(e-j-1)^{d_{j+1}-1},\\
\tilde{\lambda}^{1'}=(e-2)^{d_0+1}(e-1)^{d_1}\cdots(e-j+1)^{d_{j-1}}(e-j)^{d_j-1}(e-j-1)^{d_{j+1}-1}.
\end{array}\right. $

Assume now that we are in case (iii). Then $d_0$ is even. Let $d=d_0+d_1$. We have $\tilde{\lambda}_i=\tilde{\lambda}_i'$ for all $i\geq d+2$, since $\lambda_i=\lambda_i'$ and $\chi(\lambda_i)=\chi(\lambda_i')$ for $i\geq d+1$. Let $\tilde{\lambda}^1=(\tilde{\lambda}_1,\ldots,\tilde{\lambda}_{d+1})$ and $\tilde{\lambda}^{1'}=(\tilde{\lambda}'_1,\ldots,\tilde{\lambda}'_{d+1})$. If $\chi(e-2)<f$, then $\tilde{\lambda}_{d+1}=\tilde{\lambda}_{d+1}'$ since either $\chi(\lambda_d)-\lambda_d+\chi(\lambda_{d+1})\leq 1$ and $\chi'(\lambda_d')-\lambda_d'+\chi(\lambda_{d+1}')\leq 1$ or $d$ is odd (this happens only when $e=2$).   We have

\quad in case (iii) $\left\{\begin{array}{ll}\tilde{\lambda}^1=e^{d_0}(e-1)^{d_1+1},\ \tilde{\lambda}'=(e-1)^{d_1+1}(e-2)^{d_0}&\text{ if }\chi(e-2)=f\\\tilde{\lambda}^1=e^{d_0}(e-1)^{d_1}\tilde{\lambda}_{d+1},\ \tilde{\lambda}'=(e-1)^{d_1}(e-2)^{d_0}\tilde{\lambda}_{d+1}&\text{if }\chi(e-2)<f\end{array}\right..$

Now it is easy to see that \ref{pf-char2} (b) holds for $\rc\in\Omega_\Lg^2$ (we use also \ref{pf-char2-2} (b)). This completes the proof of Proposition \ref{prop-psi}.

\section{special pieces}\label{sec-spp}
We say that a unipotent class (resp.  nilpotent orbit) $\rc$ is special if $\gamma_G^p(\rc)$ (resp. $\gamma_\Lg^p(\rc)$) is a special character of $\bW$ (see \cite{Lu1,Lu7}). If $G$ is of type $B_n$ or $C_n$, then $(\alpha,\beta)\in\cP_2(n)$ is special if and only if $\alpha_{i+1}\leq\beta_i\leq\alpha_i+1$ for all $i\geq 1$; if $G$ is of type $D_n$, then $(\alpha,\beta)\in\bW^\wedge$ is special if and only if $\alpha_{i+1}-1\leq\beta_i\leq\alpha_i$ for all $i\geq 1$, in particular each degenerate character is special (see \cite{Lu1}).

Let $\rc$ be a special unipotent class (resp. nilpotent orbit) in $G$ (resp. $\Lg$). We define the corresponding special piece $\cS_{\rc}$ to be the subset of $\cU_G$ (resp. $\cN_\Lg$) consisting of all elements in the closure of $\rc$ which are not in the closure of any special unipotent class (resp. nilpotent orbit)  $\rc'< \rc$ (see \cite{Lu7} when $p=1$). We show that a special piece is a union of unipotent (resp. nilpotent) pieces (for unipotent case, see also \cite{Lu8}). Hence $\cU_G$ (resp. $\cN_\Lg$) is partitioned into special pieces $\cS_\rc$ indexed by special unipotent (resp. nilpotent) classes $\rc$ (when $p\neq 2$, this follows from \cite{Lu7}). In the remainder of this subsection assume $p=2$.

Let $\rc\in\Omega_G^2$ (resp. $\Omega_\Lg^2$)  be a special class and let $\cS_{\rc}$ be the corresponding special piece. There exists $\tilde{\rc}\in\Omega^1_G$ such that $\gamma^2(\rc)=\gamma^1_G(\tilde{\rc})$.  Assume the corresponding special piece $\cS_{\tilde{\rc}}$ (in the unipotent variety of the group over $\mathbb{C}$ of the same type as $G$) is a union of the special class $\tilde{\rc}:=\tilde{\rc}^0$ and non-special classes $\tilde{\rc}^1,\ldots,\tilde{\rc}^h$, where $\tilde{\rc}^i\in\Omega_G^1$, $i\in[0,h]$. We show that

\begin{center}
(a)  $\cS_{\rc}=\sqcup_{i\in[0,h]}\Sigma^{2,G}_{\tilde{\rc}^i}$ (resp. $\cS_{\rc}=\sqcup_{i\in[0,h]}\Sigma^{2,\Lg}_{\tilde{\rc}^i}$).
\end{center}

Assume $\gamma_G^1(\tilde{\rc})=(\tilde{\alpha},\tilde{\beta})$. Let $\rc^*\in\cS_{\rc}$ and assume $\gamma^2(\rc^*)=(\alpha,\beta)$. Then $(\alpha,\beta)\leq(\tilde{\alpha},\tilde{\beta})$ and for any special $(\tilde{\alpha}',\tilde{\beta}')<(\tilde{\alpha},\tilde{\beta})$, $(\alpha,\beta)\nleq(\tilde{\alpha}',\tilde{\beta}')$.  Assume $\Phi_G(\alpha,\beta)=(\tilde{\alpha}^*,\tilde{\beta}^*)$. It follows from \ref{ssec-fiber1} (b) that  $(\tilde{\alpha}^*,\tilde{\beta}^*)\nleq(\tilde{\alpha}',\tilde{\beta}')$ and from \ref{ssec-fiber1} (c) that $(\tilde{\alpha}^*,\tilde{\beta}^*)\leq(\tilde{\alpha},\tilde{\beta})$. Hence $(\tilde{\alpha}^*,\tilde{\beta}^*)=(\tilde{\alpha}^i,\tilde{\beta}^i)$ for some $i\in[0,h]$ and thus $\rc^*\in\Sigma_{\tilde{\rc}^i}^{2, G}$ (resp. $\Sigma_{\tilde{\rc}^i}^{2, \Lg}$) (note if $\rc$ is a degenerate class, then $h=0$ (see \cite{Lu1,Lu7}) and $\rc^*=\rc$).  This shows that $\cS_{\rc}\subset\sqcup_{i\in[0,h]}\Sigma^{2,G}_{\tilde{\rc}^i}$ (resp. $\cS_{\rc}\subset\sqcup_{i\in[0,h]}\Sigma^{2,\Lg}_{\tilde{\rc}^i}$).

Now if $\rc$ is a degenerate class, then $h=0$ and the r.h.s of (a) is $\{\rc\}\subset\cS_\rc$ (see \ref{ssec-fiber1} (a)). Assume $\rc$ is not a degenerate class and assume $\gamma_G^1(\tilde{\rc}^j)=(\tilde{\alpha}^j,\tilde{\beta}^j)$, $j\in[0,h]$. Assume $\rc^*\in\Sigma^{2, G}_{\tilde{\rc}^j}$ (resp. $\rc^*\in\Sigma^{2, \Lg}_{\tilde{\rc}^j}$) and $\gamma^2(\rc^*)=(\alpha^j,\beta^j)$. Then we have $\Phi_G(\alpha^j,\beta^j)=(\tilde{\alpha}^j,\tilde{\beta}^j)$ (see \ref{ssec-fiber1} (a)).
Let $\rc'<\rc$ be another special class and assume that $\gamma^2(\rc')=(\tilde{\alpha}',\tilde{\beta}')$.   We have $(\alpha^j,\beta^j)\leq(\tilde{\alpha}^j,\tilde{\beta}^j)\leq(\tilde{\alpha},\tilde{\beta})$ (see \ref{ssec-fiber1} (b))
 and $(\alpha^j,\beta^j)\nleq(\tilde{\alpha}',\tilde{\beta}')$ (see \ref{ssec-fiber1} (c)). Thus $\rc^*\in\cS_\rc$. Hence  $\Sigma_{\tilde{\rc}^i}^{2,G}\subset\cS_{\rc}$ (resp. $\Sigma_{\tilde{\rc}^i}^{2,\Lg}\subset\cS_{\rc}$), $i\in[0,h]$. The proof of (a) is completed.

It follows from (a) and the remark after \ref{ssec-comb1} (b) that 

\quad(b) {\em the number of $\tF_{p^s}$-rational points in a special piece is a polynomial in $p^s$ with integer coefficients  independent of $p$ and $s$}.

In view of the identification  $\Omega_{SO(2n+1)}^2=\Omega_{Sp(2n)}^2$ given by the special isogeny $SO(2n+1)\rightarrow Sp(2n)$, it follows from (b) that 

\quad(c) {\em the number of $\tF_{p^s}$-rational points in a special piece is a polynomial in $p^s$ that depends only on the Weyl group.}\\[10pt]
Note that (c) implies statement 6.9 (a) of \cite{Lu7} for classical groups. I am grateful to the referee for pointing out that our results imply 6.8 (a) and  6.9 (a) of \cite{Lu7} for classical groups.

\vskip 10pt {\noindent\bf\large Acknowledgement} \quad I wish to thank George Lusztig for helpful suggestions and discussions and the referee for many valuable comments. I am also grateful to Eric Sommers, Pramod N. Achar and Anthony Henderson for their interest in this work.

\end{document}